\documentclass[12pt]{article}
\usepackage{amsmath}
\usepackage{graphicx}
\newcommand{\qed}{{\unskip\nobreak\hfil\penalty50\hskip2em\vadjust{}
            \nobreak\hfil$\Box$\parfillskip=0pt\finalhyphendemerits=0\par}}

\newtheorem{thm}{Theorem}[section] 
\newtheorem{lemma}{Lemma}[section] 

\newtheorem{definition}{Definition}[section]
\newcommand{\bed}{\begin{definition}}
\newcommand{\eed}{\end{definition}}

\newcommand{\eps}{\epsilon}
\newcommand{\epshat}{\hat{\epsilon}}

\newcommand{\bitem}{\begin{itemize}}
\newcommand{\eitem}{\end{itemize}}

\newcommand{\goto}{\rightarrow}

\newcommand{\mmax}{\mathrm{max}}
\newcommand{\mmin}{\mathrm{min}}

\newcommand{\beqn}{\begin{equation}}
\newcommand{\eeqn}{\end{equation}}
\newcommand{\balign}{\begin{align}}
\newcommand{\ealign}{\end{align}}
\newcommand{\barF}{\bar{F}}

\newcommand{\so}{\sigma_0}
\newcommand{\uo}{\mu_0}
\newcommand{\s}{\sigma}

\newcommand{\that}{\hat{t}}

\newcommand{\hf}{{1\over 2}}
\usepackage{amssymb}
\begin{document}
\nocite{*}
\bibliography{apply}

\title{Estimating the Null and the Proportion of  non-Null Effects in
  Large-Scale Multiple Comparisons} 
\author{Jiashun Jin and   T.   Tony  Cai}
\author{Jiashun Jin$^1$ and   T.   Tony  Cai$^2$}
\footnotetext[1]{~ Department of Statistics, Purdue University, West
Lafayette, IN 47907.  The research of Jiashun Jin is partially
supported by NSF Grant DMS-0505423.}  
\footnotetext[2]{~ Department of Statistics, The Wharton School,
University of Pennsylvania, Philadelphia, PA 19104. The research of
Tony Cai was supported in part by NSF Grants DMS-0306576 and DMS-0604954.} 

\date{}
\maketitle

\begin{abstract}
An important issue raised by Efron \cite{Efron} in the context of
large-scale multiple comparisons is that in many
applications the usual assumption that the null distribution is known
is incorrect, and seemingly negligible differences in the null may
result in large differences in subsequent studies. This suggests that
a careful study of estimation of the null is indispensable.    

In this paper, we consider the problem  of estimating a null 
normal distribution, and a closely related problem,
estimation of the proportion of non-null effects.    
We develop an approach based on the 
empirical characteristic function and Fourier analysis. The estimators are shown to be  uniformly consistent over a wide class of  parameters. 
Numerical  performance of
the estimators is investigated using both simulated and real data.
In particular, we apply our procedure to the analysis of breast cancer
and HIV microarray data sets. The estimators perform favorably in
comparison to  existing methods.  
\end{abstract}

\begin{quote} 
\end{quote}
\begin{quote} \small
\end{quote}


\section{Introduction}
\label{sec:Intro}
\setcounter{equation}{0}
The analysis of massive data sets now commonly arising in scientific
investigations poses many statistical challenges not present in
smaller scale studies.  One such challenge is the need for large-scale
{\it simultaneous  testing}  or {\it multiple comparisons},
in which  thousands or even millions of hypotheses  are tested 
simultaneously. In this setting, one considers a large number of
null hypotheses $H_1, H_2, \ldots, H_n$,
and is interested in determining  which hypotheses are true and which
are not.   Associated with each hypothesis  is a test statistic.
When $H_j$ is true, the test statistic $X_j$  has  a null distribution
function (d.f.) $F_0$. That is,
 \[
(X_j  |  \mbox{$H_{j}$ is true})    \qquad \sim \qquad   F_0.  
\]
Since the pioneering work of Benjamini and Hochberg \cite{BH95}, which
introduced the False Discovery Rate (FDR)-controlling procedures,
research on large-scale 
simultaneous testing has been very active.  See, for example,
\cite{ABDJ, CJL, DJ04, Efron, Efronetal, Wasserman, Rice,  Storey2}. 

\begin{figure}[htb]
\centering
\includegraphics[height = 1.9 in, width = 2.5  in]{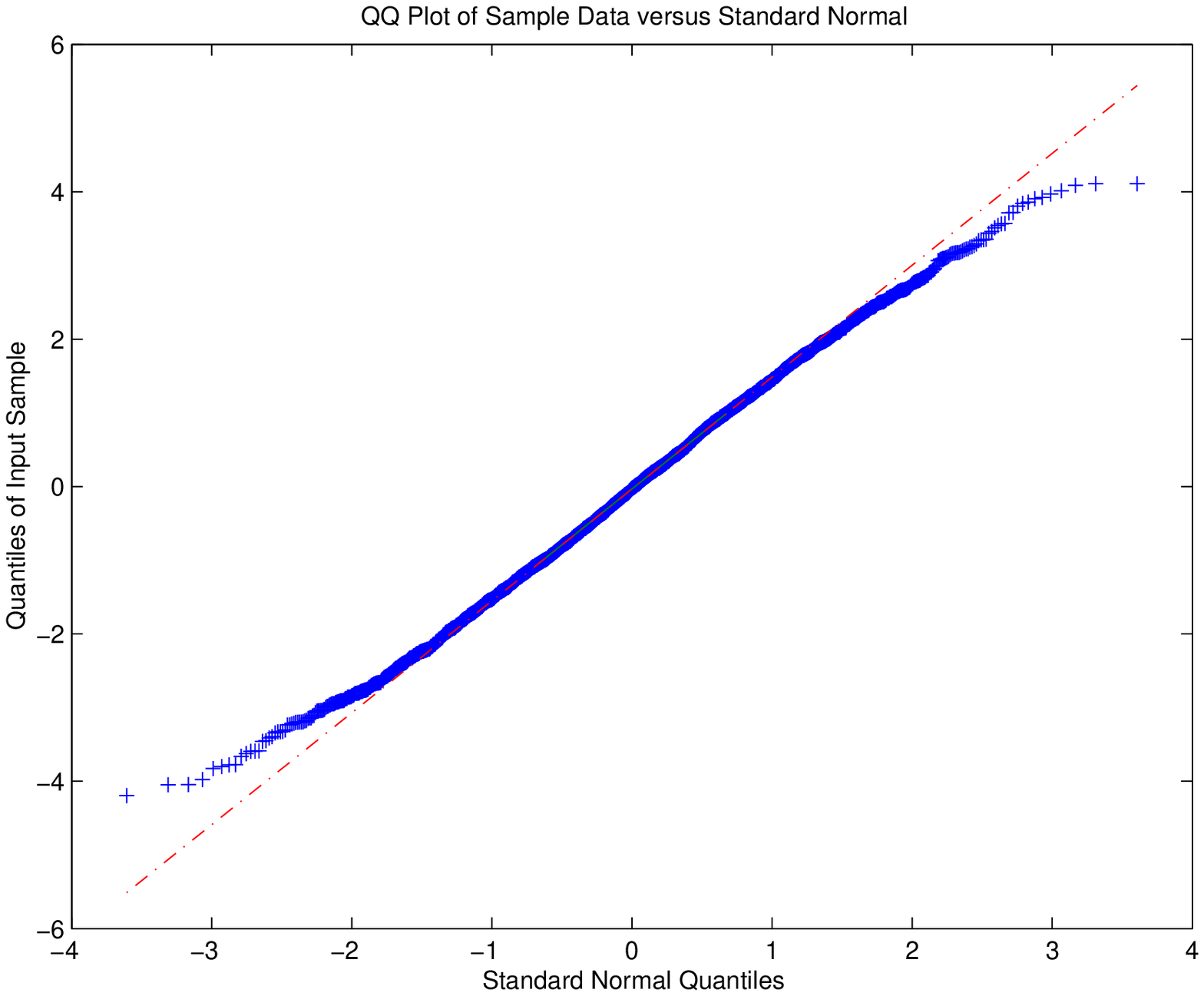}
\includegraphics[height = 1.9 in, width = 2.5  in]{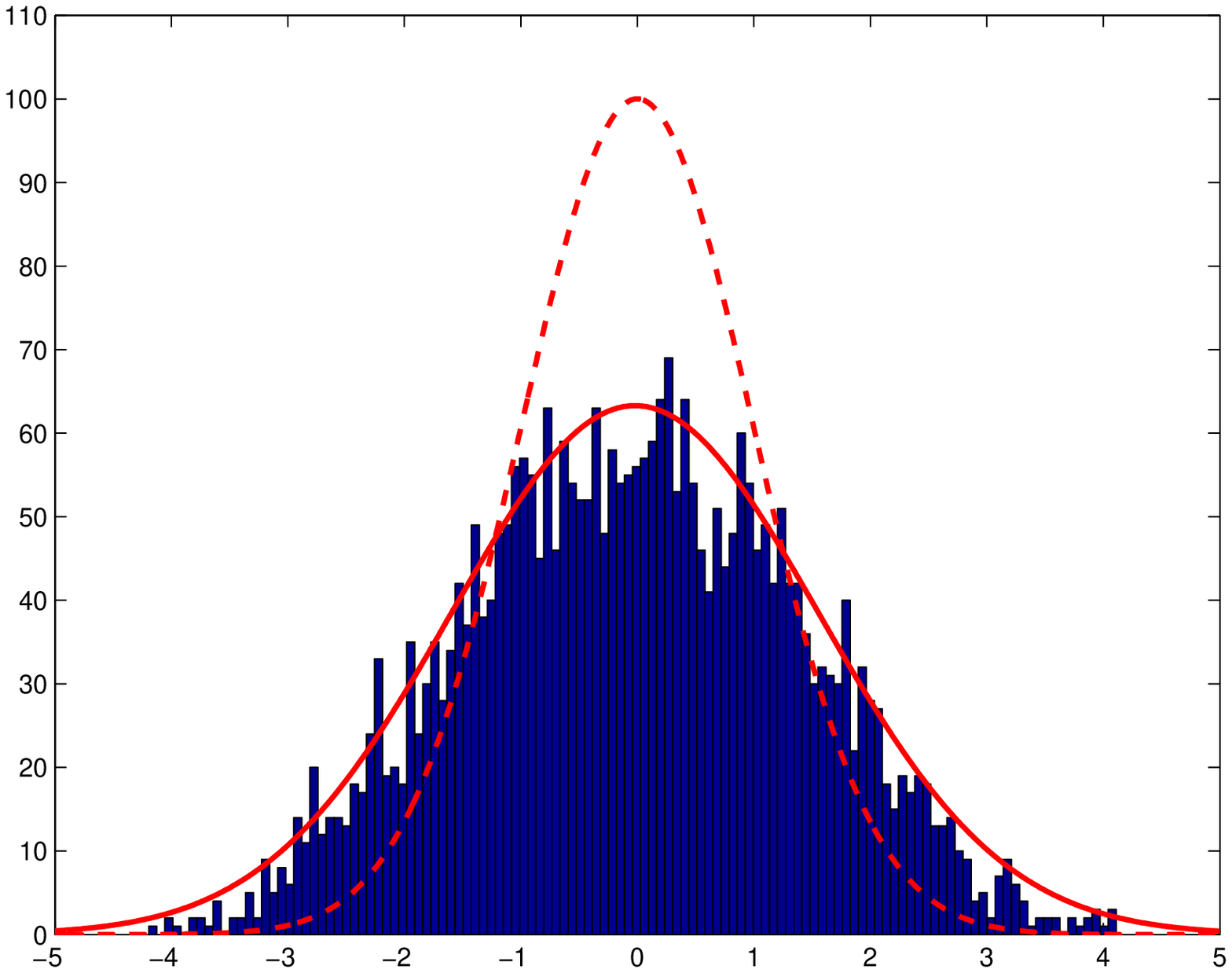}
\caption{$z$-values of microarray data on breast cancer.  Left panel:
  QQ-plot.  Right panel: histogram and density curves of  $N(0,1)$
  (dashed) and $N(-0.02, 1.58^2)$.  The plot suggests that  the null  is  
  $N(-0.02, 1.58^2)$ rather than 
  $N(0,1)$. See Efron  \cite{Efron} for further details.}
\label{fig:breastcancer}
\end{figure}

FDR  procedures are  based on the $p$-values,   which measure   the tail
probability of the null distribution. Conventionally the null
distribution is always assumed to be known. However, somewhat
surprisingly,  Efron pointed out in \cite{Efron} that in 
many applications such an assumption would be incorrect.  
Efron \cite{Efron} studied a  data set on breast cancer, in which a
gene microarray was generated for each patient in two groups, BRCA1 group and BRCA2 group. The  goal was to determine which
genes were differentially expressed between the two groups. For each gene,
a $p$-value was calculated using the classical $t$-test. For
convenience Efron chose to work on the $z$-scale through the transformation   $X_j = \bar{\Phi}^{-1}(p_j)$, 
where $\bar{\Phi} = 1 - \Phi$ is the survival function of the standard
normal distribution. Efron argued that, though theoretically the null
distribution 
should be the standard normal, empirically another null distribution
(which Efron referred to as the {\it empirical null}) is found to be more 
appropriate. In fact, he found  that  $N(-0.02, 1.58^2)$ 
is a more appropriate null than $N(0,1)$;   see Figure
\ref{fig:breastcancer}.  A similar phenomenon is also found in the
analysis of a microarray data set on HIV \cite{Efron}.  

\begin{figure}[htb]
\centering
\includegraphics[height = 1.9 in, width = 4  in]{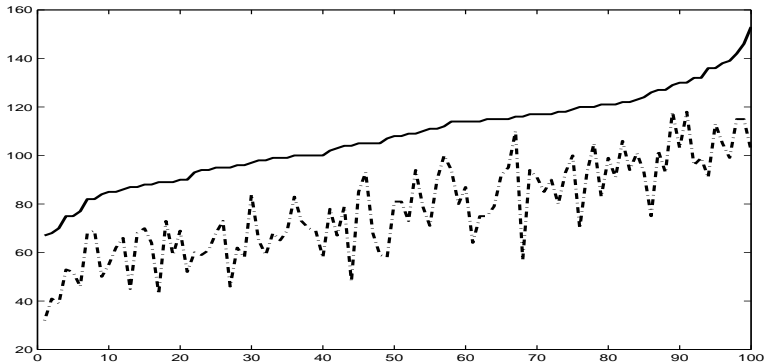}
\caption{The solid and dashed curves  represent  the 
  number of true positives for each cycle, using the true null and the
  misspecified null,  respectively.  For visualization, the numbers are sorted ascendingly with respect to those in the true null case.  }
\label{fig:wrongnull}
\end{figure}

Different choices of the null distribution can give substantially
different outcomes in simultaneous multiple testing. Even a
seemingly negligible estimation error of the null may  result in
large differences in subsequent studies.   For illustration, we
carried out an experiment which contains   $100$ independent cycles of simulations.  In  each cycle,  $9000$ samples are drawn from  $N(0, 0.95^2)$ to represent the null  effects, and  $1000$ samples are drawn from  $N(2, 0.95^2)$ to
  represent the non-null effects.     For each sample element $X_j$,    $p$-values are  calculated as  $\bar{\Phi}^{-1}(X_j/0.95)$ and
  $\bar{\Phi}^{-1}(X_j)$,   which represent the $p$-values under the {\it true}  null and the {\it misspecified} null,  respectively.        The FDR procedure is  then applied to both sets   of   $p$-values, where the FDR control parameter is set at $0.05$.    The results,  reported in Figure \ref{fig:wrongnull}, show   that  the  true positives obtained  by using $N(0, 1)$ as the null and those 
obtained  by using $N(0,0.95^2)$ as the null are considerably
different. This,
together with Efron's arguments, suggests that a careful study on
estimating the null is indispensable.    

Efron \cite{Efron} introduced a method for estimating the
null distribution based on the notion of ``sparsity.'' There are several
different ways to define sparsity \cite{ABDJ}. The most intuitive one
is that  the proportion of  non-null effects is small.   
In some applications, the case of
``asymptotically vanishing sparsity'' is of particular interest
\cite{ABDJ, DJ04}. This case refers to the situation where  the
proportion of  non-null effects tends
to zero as the  number of hypotheses grows to infinity. In such a
setting,   heuristically,    the influence of the non-null effects
becomes more and more  negligible and so the null can be reliably
estimated asymptotically.  In fact,  Efron \cite{Efron} suggested an approach which uses  the
center and half width of the central peak of the histogram for estimating the parameters of the null distribution.      

In many applications it is more appropriate to 
model the setting  as  {\it non-sparse},  i.e.,  the  proportion of 
non-null effects  does   {\it not}  tend  to zero when the number of hypotheses grows to infinity. In
such settings,  Efron's approach \cite{Efron}  does not perform
well, and  it is not hard to show that    the estimators of the null are generally
inconsistent.  Moreover, even when the setting is asymptotically 
vanishingly sparse and the estimators are  consistent,   it is
still of interest to quantify the influence of  sparsity on the
estimators, as a small error  in the null may propagate to large
errors in subsequent studies. 

Conventional methods for estimating the null parameters are based on
either moments or extreme observations  \cite{Efron, Rice, Swanepoel}.  
However, in the non-sparse case,
neither is very informative as the relevant information about  the
null is highly distorted by the non-null effects in both of them.
In this paper,   we propose a new approach  for estimating the null
parameters by using   the {\it empirical characteristic function and Fourier
analysis}  as the main tools. 
   The approach demonstrates  that the information about the null  is 
well preserved in the  high frequency Fourier coefficients,  where the distortion 
of the non-null effects is asymptotically negligible.  
The approach  integrates the   strength of several factors, including
sparsity and heteroscedasticity, and provides good estimates of the
null in a  much broader range of situations than existing approaches do.   
The resulting estimators are shown to be  uniformly consistent over a wide class of parameters  and outperform existing methods in simulations.

Beside the null distribution, the proportion of  non-null effects
is an important quantity.      For example,
the implementation of many recent
procedures requires the knowledge of both the null and the 
proportion of   non-null effects;    see \cite{Efronetal, Speed1,
   Storey2}.  Developing good estimators   for   the
proportion  is a challenging task. 
Recent work includes that of Meinshausen and Rice \cite{Rice}, Swanepoel
\cite{Swanepoel}, Cai et al.   \cite{CJL},  and Jin \cite{Jin}.    
In this paper we extend the method of Jin  \cite{Jin} to the
current setting of heteroscedasticity with an unknown null distribution.
The estimator is shown to be uniformly consistent over a  wide class of 
parameters.

In addition to the theoretical properties, numerical performance  of
the estimators is investigated using both simulated and real data.
In particular, we use our procedure to analyze the breast
cancer \cite{BreastCancer} and HIV \cite{HIV} microarray data 
that were analyzed in Efron \cite{Efron}. The results indicate
that our estimated null parameters lead to a  more reliable
identification of  differentially expressed genes than that in
\cite{Efron}.

The paper is organized as follows. In Section \ref{sec:Main}, 
after basic notations  and definitions are reviewed, the estimators of
the null parameters are defined in Section \ref{subsec:plugin}. The
theoretical properties of the estimators are investigated in Sections
\ref{sec:oracle} and \ref{subsec:Inter}. Section \ref{subsec:Depen}
discusses the extension to dependent data structures.
Section \ref{sec:Jines} treats the estimation of the proportion
of  non-null effects.   A simulation study is 
carried out in Section \ref{sec:Simul} to investigate numerical
performance. In Section \ref{sec:Appli}, we apply our
procedure to the analysis of the breast cancer \cite{BreastCancer} and HIV
\cite{HIV} microarray data. Section \ref{sec:Proof} gives proofs
of the main theorems. 

\section{Estimating the null distribution}
\label{sec:Main}
\setcounter{equation}{0}
As in Efron \cite{Efron}, we shall work on the $z$-scale and consider $n$  test 
statistics  
\begin{equation}  \label{Eqmodel2}
X_j    \sim    N(\mu_j, \s_j^2),     \qquad 1 \leq j \leq n,   
\end{equation} 
where $\mu_j$ and $\s_j$ are unknown parameters.  For a pair of  {\it null
  parameters} $\uo$ and $\so$,  
\begin{equation}   
(\mu_j,  \s_j) =  (\uo, \so) \quad \mbox{if  $H_j$ is true} 
\quad \mbox{and}\quad
(\mu_j,   \s_j)  \neq (\uo, \so) \quad \mbox{if  $H_j$ is untrue}, 
\label{Eqmodel42}
\end{equation} 
and we are interested in estimating $\uo$ and $\so$. 
We shall first consider the case in which $X_1, \ldots, X_n$ are
independent. The dependent case is considered in  Section
\ref{subsec:Depen}. 

Set $\mu = \{\mu_1, \ldots, \mu_n\}$ and $\s  =  \{\s_1, \ldots, \s_n\}$. 
Denote  the proportion  of   non-null effects by  
\begin{equation}  \label{EqDefineeps}
\eps_n    =   \eps_n(\mu,  \sigma)   =   \frac{\#\{j:\;  (\mu_j, \sigma_j) \neq (\uo, \so)\}}{n}.      
\end{equation}  
We assume $\sigma_j \geq \so$ for all $1 \leq j \leq n$. That is, 
the standard deviation of a  non-null effect   is no less
than  that of a null effect.  This is the case in a wide range of applications   \cite{Efron, Speed1}.   To make the null parameters
identifiable,  we shall assume 
\begin{equation}    \label{Eqmodel5}
\eps_n(\mu, \sigma)   \leq \eps_0,   \qquad   \mbox{for some constant $0 < \eps_0 < \hf$}.
\end{equation} 
\bed  
Fix $\eps_0  \in (0, 1/2)$,  $\uo$, and  $\so > 0$.   We say that
$(\mu, \s)$ is $(\uo, \so,\eps_0)$-eligible if (\ref{Eqmodel5}) is
satisfied and $\sigma_j \geq \so$ for all $1 \leq j \leq n$.    
\eed
Throughout this paper,  we assume  that   $(\mu, \sigma)$ is  $(\uo, \so, \eps_0)$-eligible.

\subsection{Estimating the null parameters}  
\label{subsec:plugin}

As mentioned in the Introduction, an informative approach for estimating 
the null distribution is to use the Fourier coefficients at suitable
frequencies. In the literature, Fourier coefficients  have been
frequently used for statistical inference;  see for example
\cite{Fan, Zhang}. We  now use them  to construct 
estimators for the null parameters.    
 
Introduce  the {\it empirical characteristic function}
\begin{equation}  \label{EqDefinephin}
\varphi_n(t)   = \varphi_n(t; X_1, \ldots, X_n, n)  =  \frac{1}{n} \sum_{j = 1}^n  e^{it X_j},  
\end{equation} 
and its expectation, the {\it characteristic function}
$\varphi(t)  =   \varphi(t; \mu, \sigma, n) =   \frac{1}{n}  \sum_{j =1}^n   e^{i t \mu_j  - \frac{\s_j^2 t^2}{2}}$,
where  $i  = \sqrt{-1}$.     
The characteristic function $\varphi$ naturally splits into two components,   
$\varphi(t)   =   \varphi_0(t) + \tilde{\varphi}(t)$,    where
$\varphi_0(t)   = \varphi_0(t; \mu, \sigma, n)   =  (1 - \eps_n)   \cdot e^{ i \uo  t  -\so^2 t^2/2}$ and 
\begin{equation}   \label{EqDefinephitilde}
\tilde{\varphi}(t)   =  \tilde{\varphi}(t; \mu, \sigma, n)  =  \eps_n   \cdot  \mathrm{Ave}_{\{j:\;  (\mu_j, \sigma_j) \neq (\uo, \so)\}}   \{ e^{ i \mu_j    t  -\s_j^2  t^2/2} \},      
\end{equation}
which  correspond to  the null effects and non-null effects, respectively.  
Note that the identifiability condition $\eps_n \leq \eps_0 < 1/2$ ensures that  $\varphi(t) \neq 0$ for all $t$.  
 
We now use the above functions to construct estimators for $\so^2$ and
$\uo$. For any   $t \neq 0$ and any differentiable complex-valued
function $f$ such that $|f(t)| \neq 0$,   we  define the two
functionals
\begin{equation}  \label{EqDefineso} 
\so^2(f; t)  =  - \frac{\frac{d}{dt}|f(t)|}{t \cdot  |f(t)|},\;
\uo(f; t) = \frac{\mathrm{Re}(f(t)) \cdot \mathrm{Im}(f'(t))
  - \mathrm{Re}(f'(t)) \cdot  \mathrm{Im}(f(t))} {|f(t)|^2}, 
\end{equation}
where $\mathrm{Re}(z)$ and $\mathrm{Im}(z)$ denote respectively the real and
imaginary parts of the complex number $z$.   Simple calculus
shows that evaluating the functionals at $\varphi_0$ gives the {\it
  exact}  values of $\so^2$ and $\uo$: $\so^2(\varphi_0; t)    =
\so^2$ and  $\uo(\varphi_0; t)  =   \uo$ for all $t \neq 0$.

Inspired by this,  we hope that 
for an appropriately chosen large $t$,   $\varphi_n(t) \approx \varphi(t) \approx \varphi_0(t)$,   so that  the contribution of non-null effects to the empirical 
characteristic function is negligible, which would  then give  rise to good estimates  for $\so^2$ and $\uo$.   More specifically,  we use   $\so^2(\varphi_n; t)$ and
$\uo(\varphi_n; t)$ as estimators for $\so^2$ and $\uo$, respectively, and  hope that by  choosing an  appropriate $t$,   
\begin{align}
\so^2(\varphi_n; t)   \;    \approx  \;    \so^2(\varphi; t)      \;   \approx  \so^2(\varphi_0; t)  \equiv  \so^2,      \label{EqApproxa} \\
\uo(\varphi_n; t)   \;    \approx  \;    \uo(\varphi; t)      \;   \approx  \uo(\varphi_0; t)  \equiv  \uo.        \label{EqApproxb}
\end{align} 
There is clearly a tradeoff in the choice of $t$.    
As $t$ increases from $0$ to $\infty$, the second approximations in
(\ref{EqApproxa}) and (\ref{EqApproxb}) become increasingly accurate,
but the first approximations   become more unstable because
the variances  of   $\so^2(\varphi_n; t)$ and  $\uo(\varphi_n; t)$
increase with $t$.  Intuitively, we should choose a $t$ such that 
$\varphi_n(t)/\varphi(t) \approx 1$,   
so that  $\varphi$ can be estimated with first order accuracy.  Note that  by the  central limit theorem,      
$|\varphi_n(t) - \varphi(t)|   =  O_p(\frac{1}{\sqrt{n}})$,    
so  $t$ should be chosen such that 
$\varphi(t)   \gg \frac{1}{\sqrt{n}}$.  

We introduce the  following method for choosing $t$,
which is adaptive to the magnitude of the empirical characteristic
function. For a given $\gamma \in (0, 1/2)$, set
 \begin{equation}  \label{EqDefinethatn}
 \that_n(\gamma)  =   \that_n(\gamma; \varphi_n)  =   \inf \{t: \;     |\varphi_n(t)|  = n^{-\gamma},  \; 0 \leq t \leq \log n  \}. 
\end{equation}  
Once we decide on the  frequency $t = \that_n(\gamma)$,   we have the
following  family of `plug in' estimators which are indexed by $
\gamma \in (0, 1/2)$:  
\begin{equation}
\hat{\sigma}_0^2  = \so^2(\varphi_n;  \that_n(\gamma))
\quad\mbox{and}\quad
\hat{\mu}_0   =   \uo(\varphi_n;    \that_n(\gamma)).  
\label{EqDefinesigmahat}
\end{equation}
We mention here that it will be shown later in Lemma \ref{lemma:t}
that  $\that_n(\gamma)$ is asymptotically equivalent to the 
non-stochastic quantity
\begin{equation}  \label{EqDefinetn}
t_n(\gamma)  =   t_n(\gamma; \varphi)  =   \inf \{t:\;     |\varphi(t)| = n^{-\gamma},  \;   0 \leq t \leq  \log n\},    
\end{equation}  
and that  the stochastic fluctuation of $\that_n(\gamma)$ is algebraically small
and its effect is generally  negligible.  We notice  here that  by elementary calculus,   \begin{equation}  \label{EqDefinetnadd} 
t_n(\gamma,  \varphi)   =   [\sqrt{2 \gamma \log n}/\sigma_0]  \cdot (1 + o(1)),  \qquad  n \goto \infty,  
\end{equation} 
where $o(1)$ tends to $0$    uniformly for all $\varphi$ under consideration.  


\subsection{Uniform consistency of the estimators}
\label{sec:oracle}

We now show that the estimators $\hat{\sigma}_0^2$ and $\hat{\mu}_0$ given in  
(\ref{EqDefinesigmahat}) are consistent uniformly over a wide class of
parameters.   
Introduce two non-stochastic bridging quantities,     $\so^2(\varphi;
t_n(\gamma))$  and   $\uo(\varphi; t_n(\gamma))$,  which correspond to
$\so^2$ and $\uo$, respectively.   For each estimator,     the
estimation error can  be decomposed into two 
components: one is the stochastic fluctuation and  the other is the difference between the true parameter and its corresponding bridging quantity,
\begin{align} 
|\so^2(\varphi_n;  \that_n(\gamma))  - \so^2|  & \leq   |\so^2(\varphi_n; \that_n(\gamma)) - \so^2(\varphi; t_n(\gamma))|  +  |\so^2(\varphi; t_n(\gamma)) - \so^2|,      \label{EqApprox2a} \\
|\uo(\varphi_n;  \that_n(\gamma))  - \uo|  &\leq  |\uo(\varphi_n; \that_n(\gamma)) - \uo(\varphi; t_n(\gamma))|  +  |\uo(\varphi; t_n(\gamma)) - \uo|. \label{EqApprox2b}  
\end{align}  
We shall consider the behavior of the two components separately. 
Fix  constants $q > 0$ and $A > 0$,  and  introduce the set of  parameters 
\begin{equation}  \label{EqDefineLambda}
\Lambda_n(q,  A; \uo,  \so,  \eps_0)  =  \{\mbox{$(\mu, \sigma)$ is $(\uo, \so,\eps_0)$-eligible},   \;    M_n^{(q)}(\mu, \sigma) \leq A^q \},  
\end{equation}
where   $M_n^{(q)}(\mu, \sigma)   = \mathrm{Ave}_{\{j:\; (\mu_j, \sigma_j) \neq (\uo, \so)\}}
\{(|\mu_j-\mu_0| +  |\sigma_j^2 - \so^2|^{1/2})^q\}$.   For a constant $r$,  we say that  a sequence $\{a_n\}_{n = 1}^{\infty}$ is  $\bar{o}(n^{-r})$ if for any $\delta  > 0$,     $n^{r  - \delta} |a_n| \goto 0$ as $n \goto \infty$. 
The following theorem elaborates  the magnitude of  the stochastic  component.

\begin{thm} \label{thm:Main1}
Fix constants $\gamma, \eps_0 \in (0,1/2)$,    $q \geq 3$, and $A > 0$.   
As $n \goto \infty$,  except for an event with probability $\bar{o}(n^{-c_1})$, \[
\sup_{\{\Lambda_n(q, A;  \uo, \so, \eps_0)\}} |\so^2(\varphi_n; \that_n(\gamma))   - \so^2(\varphi; t_n(\gamma))|  \leq   3 c_2    \cdot \log^{1/2}(n)  \cdot  n^{\gamma -1/2},  
\]
\[
\sup_{\{\Lambda_n(q, A;  \uo, \so, \eps_0)\}} |\uo(\varphi_n;  \that_n(\gamma))  - \uo(\varphi;  t_n(\gamma))|  \leq  \sqrt{2 \gamma} c_2    \cdot  \log (n) 
\cdot n^{\gamma - 1/2},
\]
where $c_2 = c_2(\so, q, \gamma)  =   2 \so^2 \cdot \sqrt{\mmax\{3,  q - 1 - 2 \gamma \}}$, and 
\begin{equation}    \label{EqDefinec}
c_1 = c_1(q, \gamma)  =  \left\{
\begin{array}{ll}
(q/2 - 1 - \gamma)/2,    &\qquad   q < 4, \\
(q/2 - 1 - \gamma), &\qquad  4 \leq  q \leq 4 + 2 \gamma, \\
(q - 1 - 2 \gamma)/3,   &\qquad q > 4 + 2 \gamma.
\end{array}
\right.
\end{equation}  
\end{thm}  
Theorem \ref{thm:Main1} says that the stochastic  components    in (\ref{EqApprox2a})  and (\ref{EqApprox2b}) 
are both algebraically small,     uniformly over $\Lambda_n$.

We now consider the non-stochastic components in   (\ref{EqApprox2a}) and (\ref{EqApprox2b}).  As  defined   in (\ref{EqDefinephitilde}),    $\tilde{\varphi}(t)$ naturally factors into 
$\tilde{\varphi}(t) =  e^{i   \uo t - \so^2 t^2/2} \cdot  \psi(t)$,   
where 
 \begin{equation}  \label{EqDefinepsi}
\psi(t) = \psi(t; \mu, \sigma, n)   =    \eps_n \cdot  \mathrm{Ave}_{\{j:\; (\mu_j, \sigma_j) \neq (\uo, \so) \}}  e^{i (\mu_j  - \uo)  t   - (\s_j^2 - \so^2) t^2/2}.  
\end{equation} 
Lemma \ref{lemma:g} in Section \ref{sec:Proof} tells us that 
  there is a constant $C  > 0$ such that uniformly for all $(\uo, \so, \eps_0)$-eligible parameters $(\mu, \sigma)$, 
$|\so^2(\varphi; t_n(\gamma)) - \so^2|  \leq  C \cdot
|\psi'(t_n(\gamma))|/t_n(\gamma)$ and
$|\uo(\varphi; t_n(\gamma)) - \uo|  \leq  C \cdot   |\psi'(t_n(\gamma))|$; see details therein.    Combining these  with Theorem \ref{thm:Main1} gives the following
theorem, which is proved in  Section \ref{sec:Proof}. 

\begin{thm}  \label{thm:Main2b}
Fix constants  $\gamma,  \eps_0  \in (0,1/2)$,   $q  \geq 3$, and $A > 0$. For all $t$,     $\sup_{\{\Lambda_n(q,A; \uo,\so, \eps_0)\}}  |\psi'(t)|  \leq A \cdot \eps_0$. Moreover,  there is a constant $C = C(\gamma,   q, A, \eps_0,  \uo, \so)$ such that,  except for an event with algebraically small
  probability, for any $(\mu, \sigma) \in \Lambda_n(q,A; \uo, \so,
  \eps_0)$ and all sufficiently large $n$, 
\[
|\so^2(\varphi_n;  \that_n(\gamma)) - \so^2|   \;\;   \leq  \;\;  C  \biggl(  \frac{|\psi'(t_n(\gamma))|}{\sqrt{\log n}}   +   \log^{1/2}(n) \cdot n^{\gamma -1/2} \biggr), 
\]
\[
 |\uo(\varphi_n;  \that_n(\gamma)) - \uo|   \;\;   \leq \;\;   C \biggl(|\psi'(t_n(\gamma))|  +  \log(n) \cdot n^{\gamma -1/2}\biggr).
\]
Consequently, $\so^2(\varphi_n; \that_n(\gamma))$  is   uniformly 
consistent for $\so^2$   over $\Lambda_n(q,A; \uo, \so, \eps_0)$.  Additionally,   if  $\psi'(t_n(\gamma)) = o(1)$, then  $\uo(\varphi_n; \that_n(\gamma))$ is  consistent for $\uo$ as well. 
\end{thm}
We remark here that $\uo(\varphi_n; \that_n(\gamma))$ is  uniformly
consistent for $\uo$  over any subset  $\Lambda^*_n \subset \Lambda_n$
with   $\sup_{\{\Lambda^*_n\}}\{|\psi'(t_n(\gamma))| \}=o(1)$. 
Although at first glance the convergence rates  are
relatively slow, they are in fact much faster in many situations.
 
\subsection{Convergence rate:   examples and discussions}
\label{subsec:Inter} 

We now show that under mild conditions the  convergence
rates of $\so^2(\varphi_n; \that_n(\gamma))$ and $\uo(\varphi_n;
\that_n(\gamma))$ can be significantly improved, and sometimes are
algebraically fast.

\medskip\noindent
{\bf Example I}.  {\it Asymptotically vanishing sparsity}.
Sparsity is a natural phenomenon found in many scientific fields such
as genomics, astronomy, and  image processing.  As mentioned
before,   asymptotically vanishing sparsity refers to the case
where  $\eps_n(\mu, \sigma) \goto 0$ (as $n \goto \infty$).  Several  models for sparsity have been considered in the literature, and among them are  {\it moderately  sparse} and  {\it very sparse},  where  $\eps_n = n^{-\beta}$ for some parameter $\beta$ satisfying   $\beta \in (0, 1/2)$ and $ \beta \in (1/2,1)$, respectively \cite{ABDJ,DJ04}. Lemma \ref{lemma:g} shows that uniformly over $\Lambda_n$,    
$|\psi'(t_n(\gamma))|     \leq        O(\eps_n(\mu, \sigma))$.  
Theorem \ref{thm:Main2b} then yields the fact that the estimation errors of
$\so^2(\varphi_n; \that_n(\gamma))$  and $\uo(\varphi_n;
\that_n(\gamma))$  are   algebraically small for  both the  moderately sparse case and the  very sparse case. 

\medskip\noindent
{\bf  Example II}. {\it  Heteroscedasticity}.   It is natural
in many applications to find that a non-null effect has an elevated
variance. A test statistic consists of  two components, signal and
noise.    An elevation of variance occurs  when the signal component
contributes extra variance. Denote the minimum elevation of the
variance for the non-null effects by 
\begin{equation}  \label{EqDefinetaun}
\tau_n   =  \tau_n(\mu, \sigma) =    \mmin_{\{j: \;  (\mu_j, \sigma_j) \neq (\uo, \so) \}}  \{\s_j^2 - \so^2\}.   \end{equation} 
Lemma \ref{lemma:g}  shows that       
$|\psi'(t_n(\gamma))|  \leq  O(\eps_n   e^{- \gamma  \log(n)   \tau_n(\mu, \sigma)})$.   
So  $\psi'(t_n(\gamma))  = o(1)$  if, say,   $\tau_n \geq  \frac{\log \log n}{\log n}$, 
and   $\psi'(t_n(\gamma))$ is algebraically small  if $\tau_n \geq c_0$ for some constant $c_0  > 0$.

\medskip\noindent
{\bf  Example III}.   {\it  Gaussian hierarchical model}.      The
Gaussian hierarchical model is  widely used  in statistical inference,  
as well as  in microarray analysis;    see
Efron \cite{Efron}, for example.    A  simple  version of the model is  where 
 $\sigma_j  \equiv \so$   and the  means $\mu_j$ associated with non-null effects  are  modeled as samples from  a density function $h$,    $(\mu_j |\mbox{$H_j$ is untrue})   \stackrel{iid}{\sim}   h$.  
It is not hard to show   that $|\psi'(t_n(\gamma))|  \leq   \eps_n \cdot   |\int e^{i t_n(\gamma) u} [(u - \uo) h(u)] du |$,   
where the integral is  the Fourier transform of the function
$(u-\uo)h(u)$ at frequency $t_n(\gamma)$. By the
Riemann-Lebesgue Lemma \cite{Mallat},   $|\psi'(t_n(\gamma))|   =
o(t_n^{-k}(\gamma))$ if the $k$-th derivative  of $h(u)$   is
absolutely integrable. In particular, if   $h$ is  Gaussian, say $N(a,b^2)$,  then $|\psi'(t_n(\gamma))|   \leq  O(\eps_n \cdot  |t_n(\gamma)|   \cdot  n^{- \gamma  b^2})$ and  is algebraically small.

We note here that  sparsity, heteroscedasticity,   
and the smoothness of $h$ can occur at the same time, which makes
the convergence  even faster.    
In a sense,  our approach combines the strengths of sparsity,
heteroscedasticity,  and the smoothness of the density $h$.   The approach can
thus be viewed as an extension of Efron's approach, as it  is
consistent not only in  the asymptotically vanishingly sparse case, but also in  many 
interesting non-sparse cases.    Additionally,  in  the asymptotically
vanishingly sparse case,  the convergence  rates of our estimators
can be substantially faster than those of Efron.   For example,    this
may occur when the data set is both  sparse  and    heteroscedastic.     

\medskip\noindent
{\bf Remark:} The theory developed in Sections  \ref{subsec:plugin} -  \ref{subsec:Inter}   can be naturally extended to the
Gaussian hierarchical model, which is the Bayesian counterpart of
Model (\ref{Eqmodel2})-(\ref{Eqmodel42})  and has been widely used in
the literature;    see for example \cite{Efron, Wasserman}.    The
model treats  the test statistics $X_j$  as samples  from a two-component Gaussian mixture:    
\begin{equation}  \label{EqBayesmodel2}
X_j   \sim     (1 - \eps) N(\uo, \so^2) + \eps N(\mu_j, \s_j^2),    \qquad 1 \leq j \leq n,   
\end{equation} 
where  $(\mu_j, \s_j)$ are  samples from a bivariate 
distribution $F(\mu, \sigma)$.    The previous results can be
naturally extended to this model.

\subsection{Extension to dependent data structures}   
\label{subsec:Depen} 

We now consider the proposed approach for dependent data. As the discussions are similar, we focus on $\so^2(\varphi; \that_n(\gamma))$. Recall that the estimation error splits into a stochastic component and a non-stochastic component, 
$|\so^2(\varphi_n; \that_n(\gamma))  - \so^2|   \leq  |\so^2(\varphi_n; \that_n(\gamma)) - \so^2(\varphi; t_n(\gamma))|  + |\so^2(\varphi; t_n(\gamma))  - \so^2|$.  
Note that the non-stochastic  component only contains marginal effects
and  is unrelated to  dependence structures.  We thus need only  to  
study the stochastic component, or to extend Theorem \ref{thm:Main1}.
In fact, once  Theorem \ref{thm:Main1} is extended to the dependent
case, the extension of Theorem \ref{thm:Main2b} follows directly by
arguments similar to those given in the proof of Theorem \ref{thm:Main2b}.  
For reasons of space, we shall focus on two dependent structures: the
strongly ($\alpha$)-mixing case and  the short-range dependent case.     Denote the strongly mixing coefficients by  $\alpha(k) =  \sup_{\{1 \leq t \leq n\}}  \alpha(\sigma(X_s,     s\leq t), \sigma(X_s, s \geq t + k))$, where  $\sigma(\cdot)$ is the $\sigma$-algebra generated by the random variables specified in the  brackets, and  $\alpha(\Sigma_1, \Sigma_2) = \sup_{\{E_1 \in \Sigma_1, E_2 \in \Sigma_2\}}    |P\{E_1 \cap E_2\}  - P\{E_1\}  P\{E_2\}|$ for any two 
$\sigma$-algebras $\Sigma_1$ and $\Sigma_2$.   In the strongly mixing case,  we suppose that  $\alpha(k) \leq B k^{-d}$ for some positive constants $B$ and $d$.    
In the short-range dependent case,  we suppose $\alpha(k) = 0$ when $k
\geq n^{\tau}$ for some constant $\tau \in (0,1)$.   

Now, fix constants $a > 0$, $B > 0$,  $q \geq 3$,  and $A > 0$,   introduce the following set of parameters which we denote by  $\tilde{\Lambda}_n(a,B, q,A) = \tilde{\Lambda}_n(a,B, q,A; \eps_0, \uo, \so)$:
\[
\{(\mu, \sigma)   \in \Lambda_n(q,A; \uo, \so, \eps_0),    \;  \mmax_{\{1 \leq  j \leq n\}}\{|\mu_j| + |\sigma_j|\} \leq B \log^{a}(n)   \}. 
\] 
Note that this technical condition is not essential and can be relaxed.
The following theorem treats the strongly mixing case and is proved in \cite[Section 7]{JC}.
\begin{thm}  \label{thm:Depen} 
Fix $d > 1.5$, $q \geq 3$,    $\gamma \in (0, \frac{d-1.5}{2d + 2.5})$,    $A > 0$,  $a > 0$, and $B
> 0$. Suppose $\alpha(k) \leq B k^{-d}$  for all $1 \leq k \leq n$.     As $n \goto \infty$,   uniformly for all $(\mu, \sigma) \in \tilde{\Lambda}_n(a,B,q,A)$,  except for  an event with asymptotically vanishing probability,  
\[
|\so^2(\varphi_n; \that_n(\gamma)) - \so^2(\varphi;  t_n(\gamma))| \leq \bar{o}(n^{\gamma - 1/2}),  \;\;   |\uo(\varphi_n; \that_n(\gamma)) - \uo(\varphi;  t_n(\gamma))| \leq \bar{o}(n^{\gamma - 1/2}).
\]   
\end{thm} 
An interesting question is whether this result holds for all $\gamma \in (0,1/2)$;  we leave this for future study. 
The following theorem concerns the short-range dependent case, whose
proof is similar to that of Theorem \ref{thm:Depen} and is thus  omitted.  
\begin{thm}  
Fix $q \geq 3$,   $\tau \in (0,1)$,    $\gamma \in (0, \frac{1 - \tau}{2})$,   $A > 0$,   $a >
0$, and $B > 0$.  Suppose $\alpha(k)= 0$  for all $  k \geq n^{\tau}$.     As $n \goto \infty$,  uniformly for all $(\mu, \sigma) \in \tilde{\Lambda}_n(a,B,q,A)$,   except for an event with  asymptotically vanishing probability, 
\[
|\so^2(\varphi_n; \that_n(\gamma)) - \so^2(\varphi;  t_n(\gamma))| \leq \bar{o}(n^{\gamma - \frac{1-\tau}{2}}),  \;\;    |\uo(\varphi_n; \that_n(\gamma)) - \uo(\varphi;  t_n(\gamma))| \leq \bar{o}(n^{\gamma - \frac{1-\tau}{2}}).
\]   
\end{thm} 

We mention that consistency for  more general dependent settings is  possible provided the following two key requirements are satisfied.  First, there is an exponential type inequality for the tail probability of $|\varphi_n(t) - \varphi(t)|$ for all $t \in (0, \log n)$;  we use Hoeffding's inequality in the proof for  the  independent case,  and use \cite[Theorem 1.3]{Bosq} in the proof of Theorem \ref{thm:Depen}.  Second,     the standard deviation  of $\varphi_n(t_n(\gamma))$ has  a smaller order than that of $\varphi(t_n(\gamma))$, so that the approximation $\varphi_n(t_n(\gamma))/\varphi(t_n(\gamma)) \approx 1$ is accurate to the first order.

\section{Estimating the proportion of non-null effects} 
\label{sec:Jines}
\setcounter{equation}{0} 
The development of useful estimator for the proportion of  non-null
effects, together with the corresponding statistical analysis, poses
many challenges. Recent work includes those of Meinshausen and Rice \cite{Rice}, Swanepoel \cite{Swanepoel},
Cai, et al.   \cite{CJL}, and Jin \cite{Jin}. See also \cite{Efronetal, Wasserman}.
The first two approaches only provide consistent
estimators under a condition which Genovese and Wasserman call
``purity'' \cite{Wasserman}.   These approaches do not perform well
in the current setting  as the purity  condition is not satisfied;
see Lemma \ref{lemma:impure} for details.     Cai et al.   \cite{CJL}
largely focuses on a  very  sparse setting,
and so a more specific model is needed.  Jin \cite{Jin} 
considers estimating the proportion of nonzero normal 
means but concentrates on the  homoscedastic case with known
null parameters.  This motivates 
a careful study of estimation of the proportion in the current  setting.       

We begin by first assuming  that the null parameters are  known.  In
this case the approach of Jin \cite{Jin} can be extended to
the heteroscedastic setting here.    Fix $\gamma \in
(0,\frac{1}{2})$. The 
following estimator is proposed in \cite{Jin} for the homoscedastic case:  
\begin{equation}  \label{EqDefineepshat}
\epshat_n(\gamma)  =  \epshat_n(\gamma; X_1, \ldots, X_n,n) =  \sup_{\{0 \leq   t  \leq   \sqrt{2 \gamma \log n}\}}  \{1 -   \Omega_n(t;  X_1, \dots, X_n, n)\},     
\end{equation}  
where  $\Omega_n(t;  X_1, \dots, X_n, n)  =    \int_{-1}^{1} (1 - |\xi|)   \bigl(\mathrm{Re}(\varphi_n(t;  X_1, \ldots, X_n, n) e^{-i   \uo t  + \so^2 t^2/2}) \bigr)  d\xi$.  
This estimator   continues
to be consistent for the  current  heteroscedastic case.  Set
\[
\Theta_n(\gamma; q, A, \uo, \so,   \eps_0)  =   \{(\mu, \sigma) \in \Lambda_n(q,A; \uo, \so, \eps_0),   \;  \Delta_n  \geq \frac{\log\log n}{ \log n},  \;   \eps_n(\mu, \sigma)   \geq  n^{\gamma - \frac{1}{2}}\}, 
\] 
where  $\Delta_n  = \Delta_n(\mu, \sigma)  = \mmin_{\{j:  \:  (\mu_j,  \sigma_j) \neq  (\uo, \so)\}}  \bigl\{\mmax\{|\mu_j - \uo|^2,  |\sigma_j^2 - \so^2| \}  \bigr\}$.
\begin{thm}   \label{thm:esteps}
For any $\gamma   \in (0,1/2)$,   $q \geq 1$, and  $A > 0$,  except for an event with  algebraically  small  probability, 
$\lim_{n \goto \infty}   \big( \sup_{\{\Theta_n(\gamma; q, A, \uo, \so,  \eps_0)\}}   \{|\frac{\epshat_n(\gamma)}{\eps_n(\mu, \sigma)} -1  |  \} \big) = 0$. 
\end{thm}   
Roughly speaking,  the estimator is consistent if the proportion is asymptotically larger than $1/\sqrt{n}$.   The case where  the proportion is asymptotically smaller than $1/\sqrt{n}$ is very challenging,  and usually  it is very hard to construct consistent estimates  without a more specific model;    see \cite{CJL, DJ04} for more discussion.

We now turn to the case  where  the null parameters $(\uo, \so)$ are
unknown.   A natural approach is to first use  the proposed procedures
in Section \ref{subsec:plugin} to obtain  estimates for $\uo$ and
$\so$, say $\hat{\mu}_0$ and $\hat{\sigma}_0$,  and  then plug them
into (\ref{EqDefineepshat}) for estimation of the proportion.  This yields the estimate 
$
\epshat_n^*(\gamma;  \hat{\mu}_0,  \hat{\sigma}_0)  =  \epshat_n^*(\gamma; \hat{\mu}_0,  \hat{\sigma}_0, X_1, \ldots, X_n, n).   
$
Theorem  \ref{thm:estepsAdd} below describes how 
$(\hat{\sigma}_0, \hat{\mu}_0)$  affects the estimation accuracy of $\epshat_n^*$.
\begin{thm}   \label{thm:estepsAdd}
Fix $\eps_0 \in (0,1/2)$, $\gamma \in (0, 1/2)$,   $q \geq 1$, and  $A > 0$.   As $n \goto \infty$,  suppose that  except for an event  $B_n$ with  algebraically small  probability,  
$\mmax\{|\hat{\mu}_0 - \uo|^2,  |\hat{\sigma}_0^2 - \so^2|  \}   =
o(\frac{1}{\log n})$. Then there are a constant $C = C(\gamma, q,A,\uo, \so, \eps_0) > 0$ and an event  $D_n$
with algebraically  small probability,   such that over $B_n^c \cap D_n^c$ 
\[
|\epshat_n^*(\gamma; \hat{\mu}_0, \hat{\sigma}_0)  - \epshat_n(\gamma)|   \;\;\;  \leq     \;\;\;   C\cdot \bigl[ \log^{-3/2} (n) \cdot n^{\gamma - 1/2}   +    \log n \cdot |\hat{\sigma}_0^2 - \so^2|  +  \sqrt{\log n} \cdot  |\hat{\mu}_0 - \uo|  \bigr].  
\] 
\end{thm}  
 
Results in previous  sections show that,   under mild conditions,  the
estimation errors of $(\hat{\mu}_0, \hat{\sigma}_0)$ are
algebraically small, and  so is 
$\epshat_n^*(\gamma) -\epshat_n(\gamma)$. In the non-sparse case,  such differences are negligible and   both  $\epshat_n(\gamma)$ and 
$\epshat_n^*(\gamma)$ are consistent.    The
sparse case,  especially when the proportion
is algebraically small, is more subtle.  In this case a more
specific model is often
needed. See Cai et al.   \cite{CJL}.
 
We now compare our procedure with those in  Meinshausen and Rice
\cite{Rice} and in Cai et al.  \cite{CJL}.     
We begin by introducing the aforementioned purity condition.  
If we model the $p$-values of the test statistics as samples from a
mixing density,    $(1 - \eps) U(0,1) + \eps   h$,     where $U(0,1)$ and $h$ are the marginal densities of the $p$-values for the null effects and non-null effects respectively. 
 The purity condition is defines as 
$\mathrm{essinf}_{\{0 < p < 1\}}   h(p)  = 0$.    
Meinshausen and Rice \cite{Rice} propose a confidence lower bound for
$\epsilon$ that is valid for all $h$. Despite this advantage, the
lower bound is generally conservative and inconsistent. In fact, the
purity condition is necessary for the lower bound to 
be consistent.  Similar results can be found in Genovese and Wasserman
\cite{Wasserman}.   
Unfortunately,   the purity condition generally does not hold in our
settings. 
\begin{lemma}  \label{lemma:impure}
Let the test statistics $X_j$ be given as in (\ref{EqBayesmodel2}).
If the marginal distribution $F(\mu, \sigma)$ satisfies  either 
  $P_{F}  \{ \s > 1\} \neq 0$ or   $P_{F}\{\sigma = 1\} = 1$, but  $P_{F}\{\mu > 0\} \neq 0$
 and $P_F\{\mu < 0\}  \neq 0$,   then  the purity condition does not hold. 
 \end{lemma} 
 
Cai et al.    \cite{CJL}  consider  a very sparse setting for a
two-point Gaussian mixture model where  the proportion is modeled as
$n^{-\beta}$ with $\beta \in (\hf, 1)$.   
Their estimator is consistent  whenever consistent
estimation is possible, and it attains the optimal rate of convergence.  
 In a sense, their approach complements our method: the 
former deals with a very sparse but more specific model,   and  the latter deals with  a more general model where the level of sparsity   is much lower.

\section{Simulation  experiments}
\label{sec:Simul}
\setcounter{equation}{0}

We now turn to the numerical performance of our
estimators of the null parameters. The goal for the
simulation study is three-fold: to investigate how different choices of
$\gamma$ affect the estimation errors,    to compare the performance
of our approach with that in Efron \cite{Efron}, and to investigate the  performance of the proposed approach for dependent data.    We leave the study for real data to Section \ref{sec:Appli}. 

We first investigate the effect of $\gamma$  
on the estimation errors. 
Set $\so
= 1/\sqrt{2}$ and $\uo = -1/2$ throughout this section.  
We take $n  = 10000$,   $\eps = 0.1$, and   $a =  0.75$,  $1.00$,   $1.25$, and $1.50$   for the following simulation experiment: 
\begin{description} 
\item[Step 1.]    ({\it Main Step}).     For each $a$,  first generate
  $n \eps$  pairs of  $(\mu_j, \sigma_j)$ with $\mu_j$ from $N(0,1)$
  and $\sigma_j$ from the uniform distribution $U(a, a + 0.5)$, and  then
  generate a sample from $N(\mu_j, \sigma_j^2)$ for each pair of
  $(\mu_j,\sigma_j)$.    These $n \eps$ samples represent the non-null
  effects.  In addition, generate $n \cdot (1 - \eps)$ samples from $N(\uo, \so^2)$ to represent the null effects.    
\item[Step 2.]   For the samples obtained in Step 1,  implement   
  $\hat{\sigma}(\gamma) = \so(\varphi_n; \that_n(\gamma))$ and   $\hat{\mu}_0(\gamma)  = \uo(\varphi_n;
  \that_n(\gamma))$ for each $\gamma = 0.01, 0.02, \ldots, 0.5$.    
\item[Step 3.]    Repeat Steps 1 and 2 for $100$ independent cycles.   
\end{description}  
\begin{figure}[htb]
\centering
\includegraphics[height = 1.33 in, width = 2.5  in]{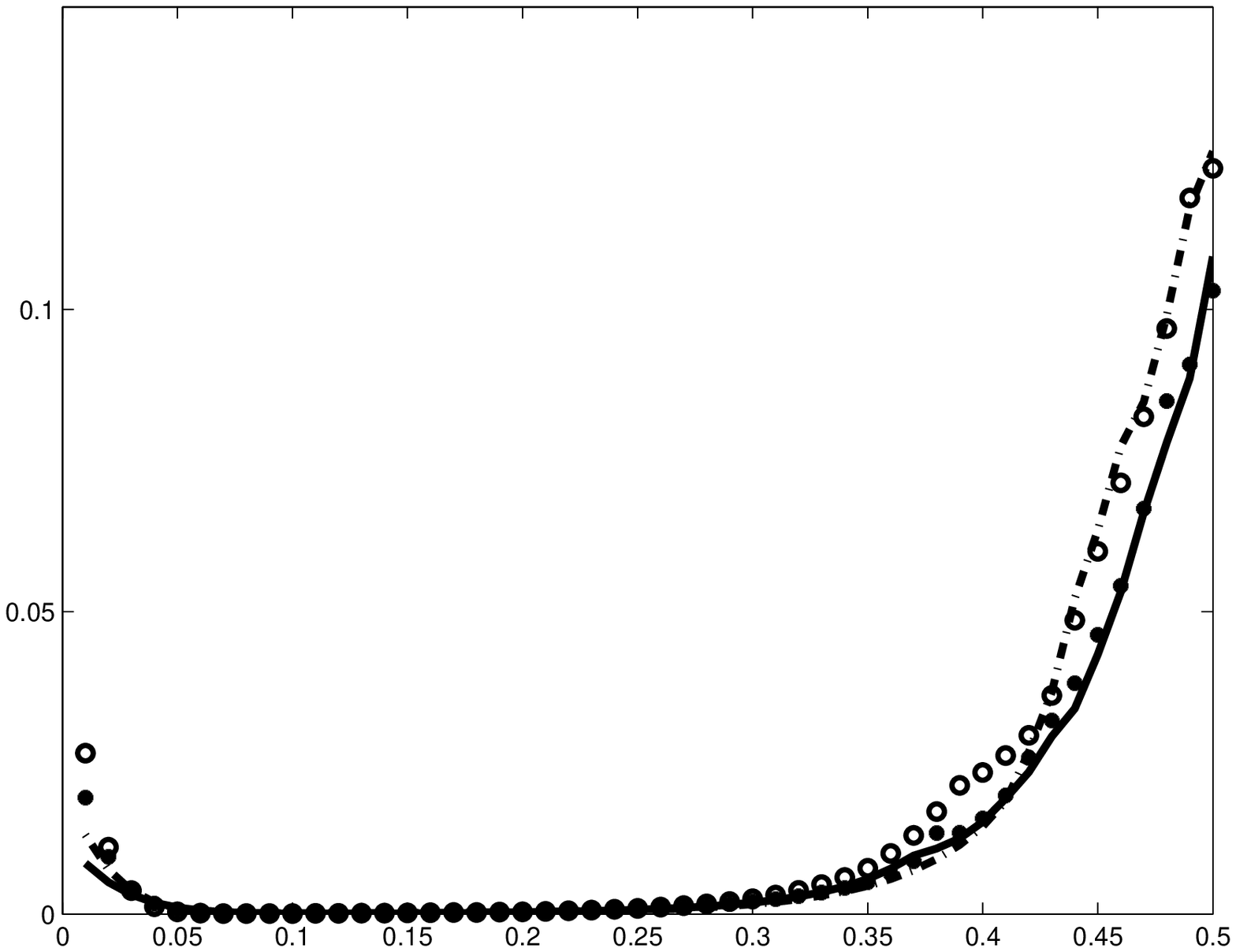}
\includegraphics[height = 1.33 in, width = 2.5  in]{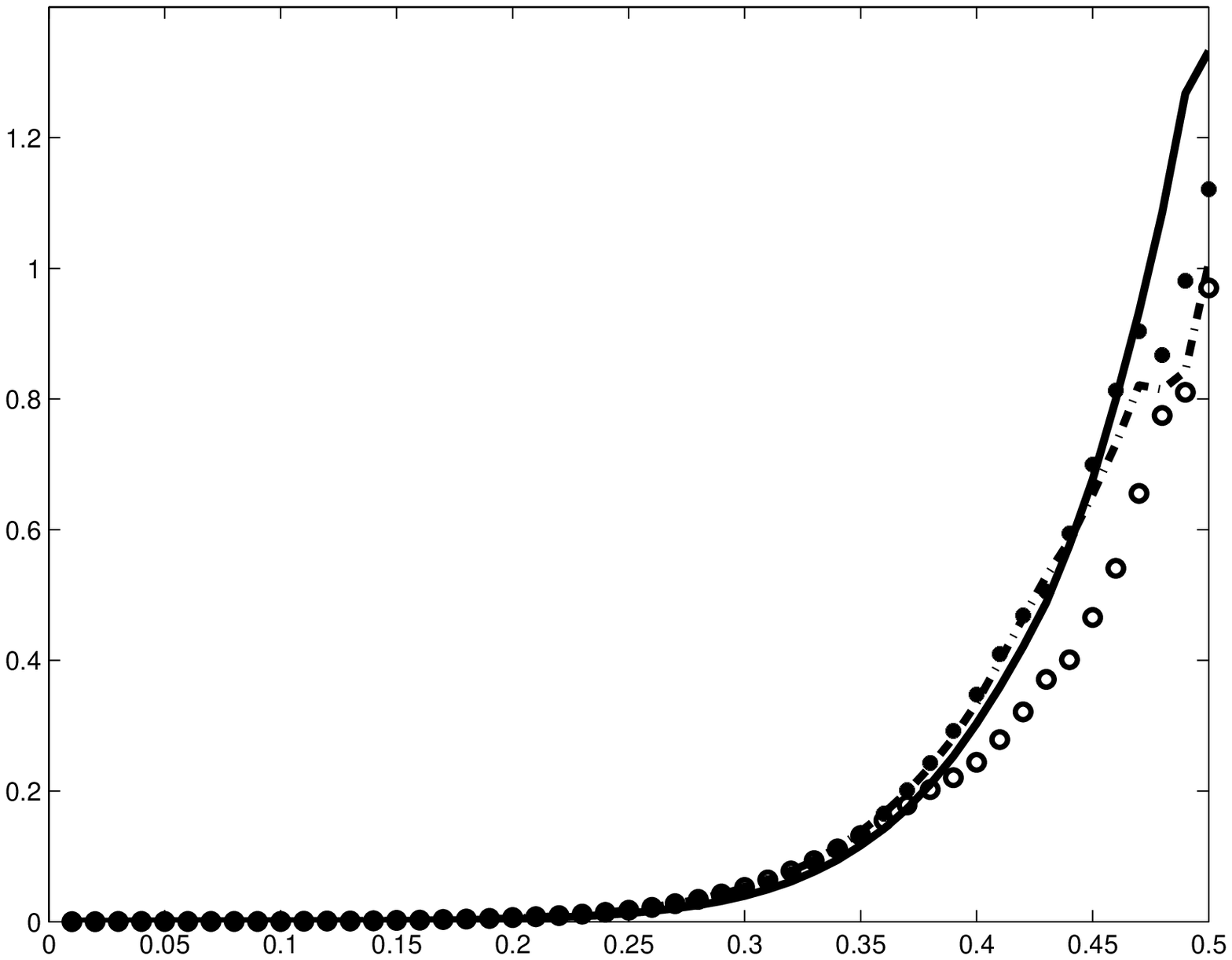}
\includegraphics[height = 1.33 in, width = 2.5  in]{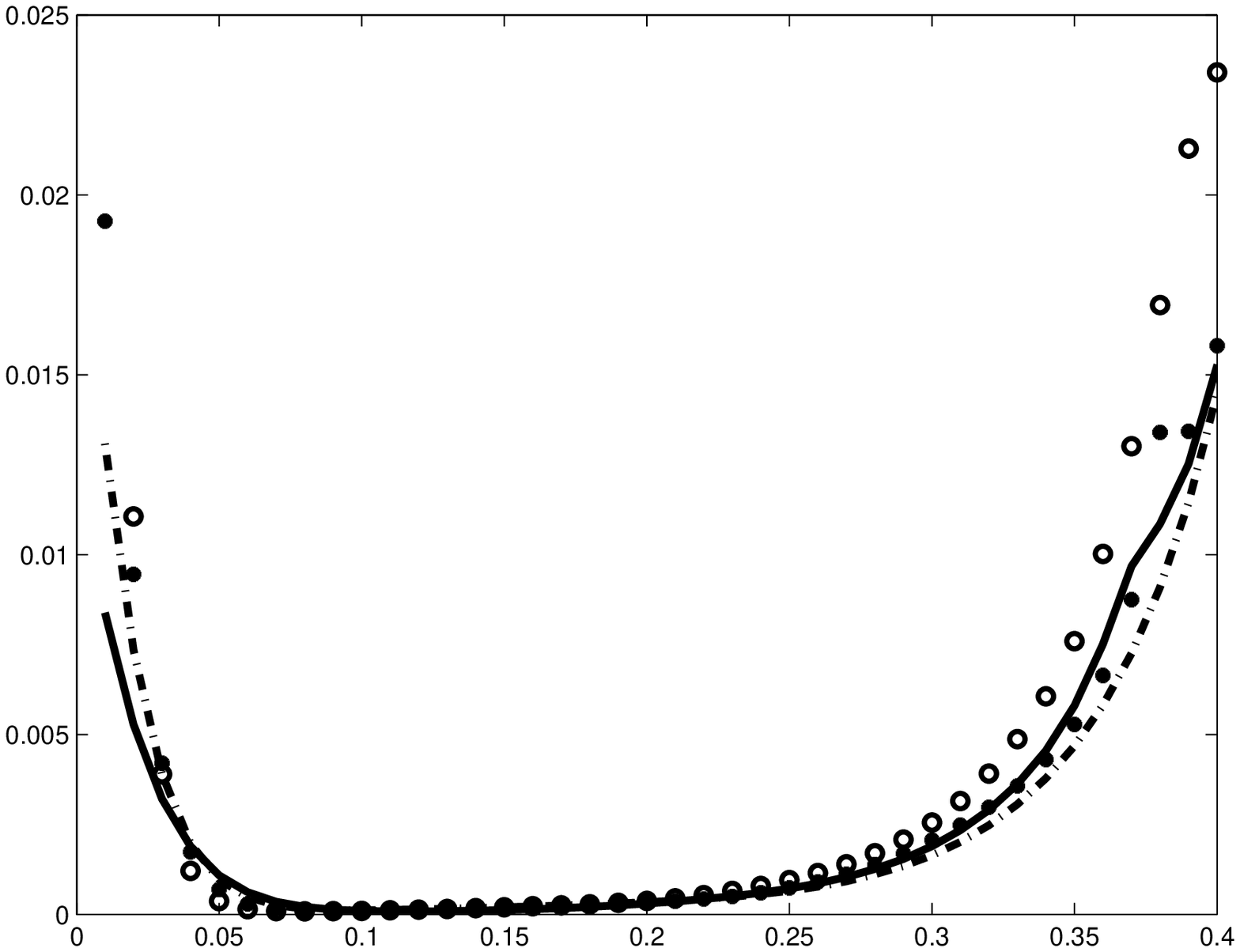}
\includegraphics[height = 1.33 in, width = 2.5  in]{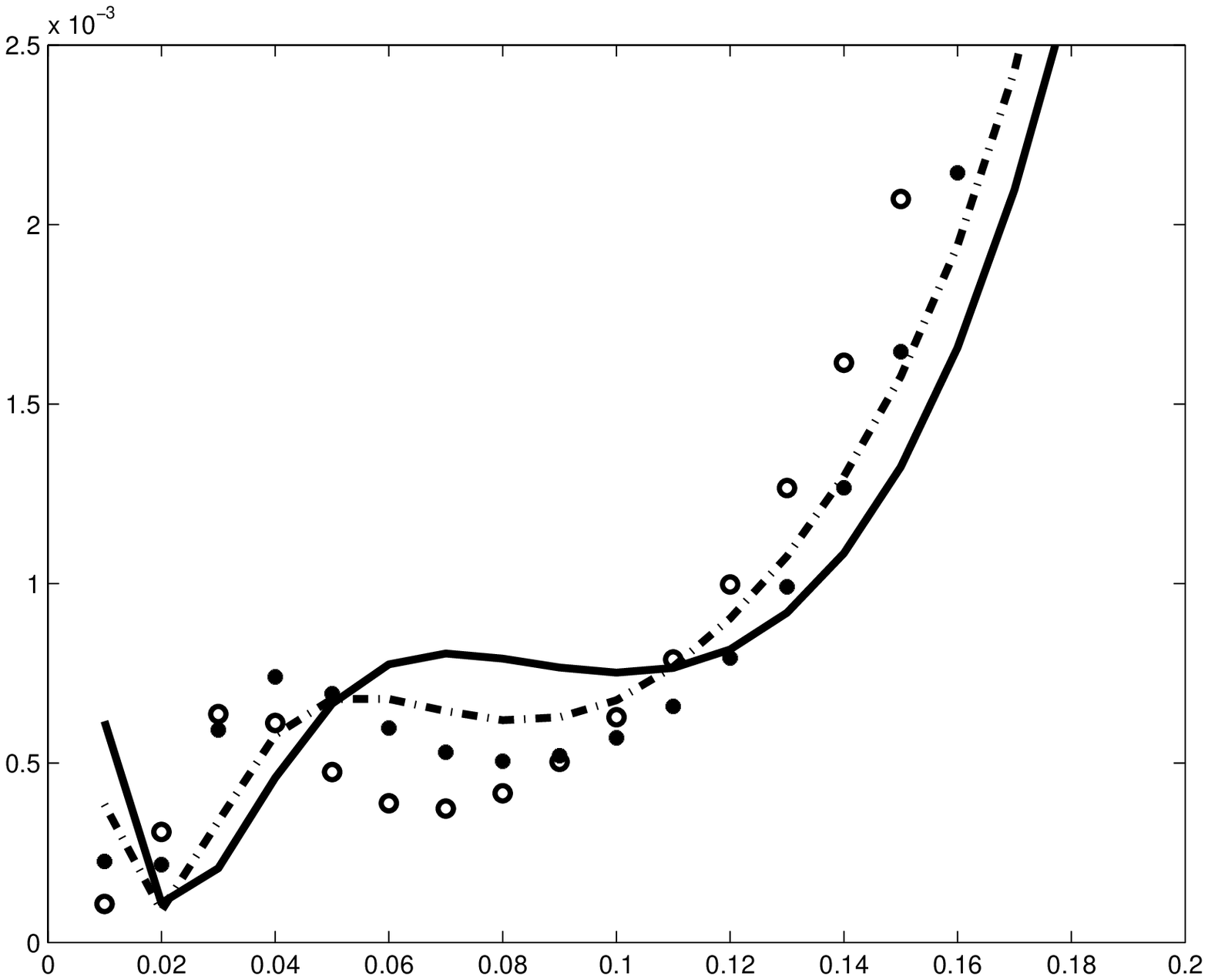}
\caption{  $x$-axis:     $\gamma$. $y$-axis:    mean squared error
  (MSE).  Top row:  MSE  for $\hat{\sigma}_0(\gamma)$ (left) and
  $\hat{\mu}_0(\gamma)$ (right).   The four different curves (solid, dashed, dot, and circle) correspond to $a = 0.75, 1.00, 1.25$, and $1.50$.     Bottom row: zoom in.}
\label{fig:pickgamma}
\end{figure}

The results, reported in Figure \ref{fig:pickgamma}, suggest that the best choice of $\gamma$ for both $\hat{\sigma}_0(\gamma)$ and $\hat{\mu}_0(\gamma)$ are  in  the range $(0.1, 0.15)$.    With $\gamma$ in this range,   the performance  of the  estimators is not very sensitive to different choices of $\gamma$,    and both estimators are  accurate.  Taking $\gamma = 0.1$, for example,   the mean squared errors for $\hat{\sigma}_0(\gamma)$ and  $\hat{\mu}_0(\gamma)$ are of magnitude $10^{-4}$ and $10^{-3}$,  respectively.     These suggest the use of the following estimators for simplicity,  where we take  $\gamma =  0.1$:
\begin{equation}  \label{EqDefinemusigmahat*}
\hat{\sigma}_0^*   = \sigma_0(\varphi_n;  \that_n(0.1)),  \qquad   \hat{\mu}_0^*    =  \uo(\varphi_n;  \that_n(0.1)).  \end{equation}

\begin{figure}[htb]
\centering
\includegraphics[height = 1.9 in, width = 4  in]{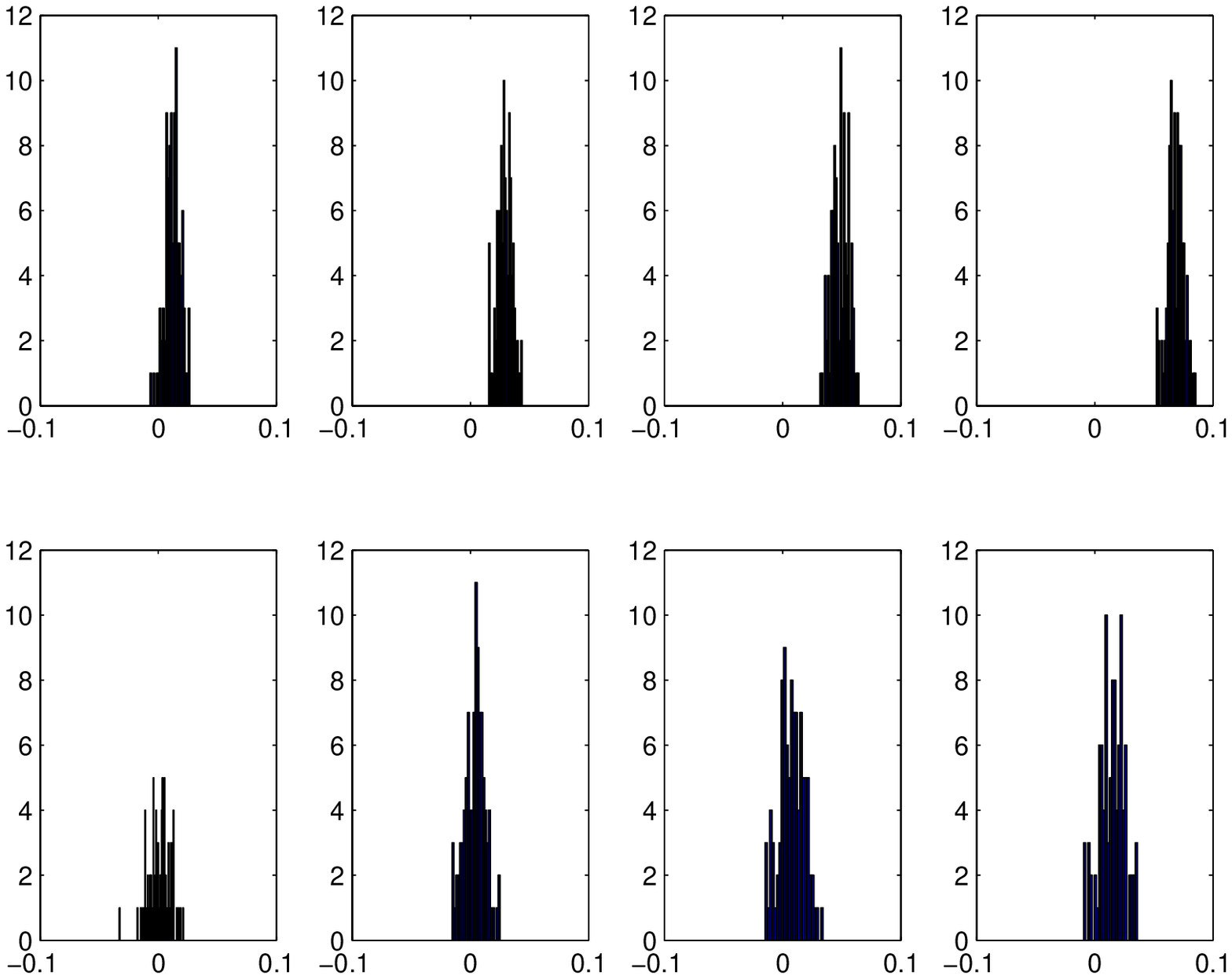}
\caption{Histograms for the  estimation errors of Efron's estimator for $\so$ (top row) and   $\hat{\sigma}_0^*$
  (bottom  row). From left to right:  $\eps =  0.05$, $0.10$, $0.15$, and $0.20$.}
\label{fig:sCompare}
\end{figure}

\begin{figure}[htb]
\centering
\includegraphics[height = 1.9 in, width = 4  in]{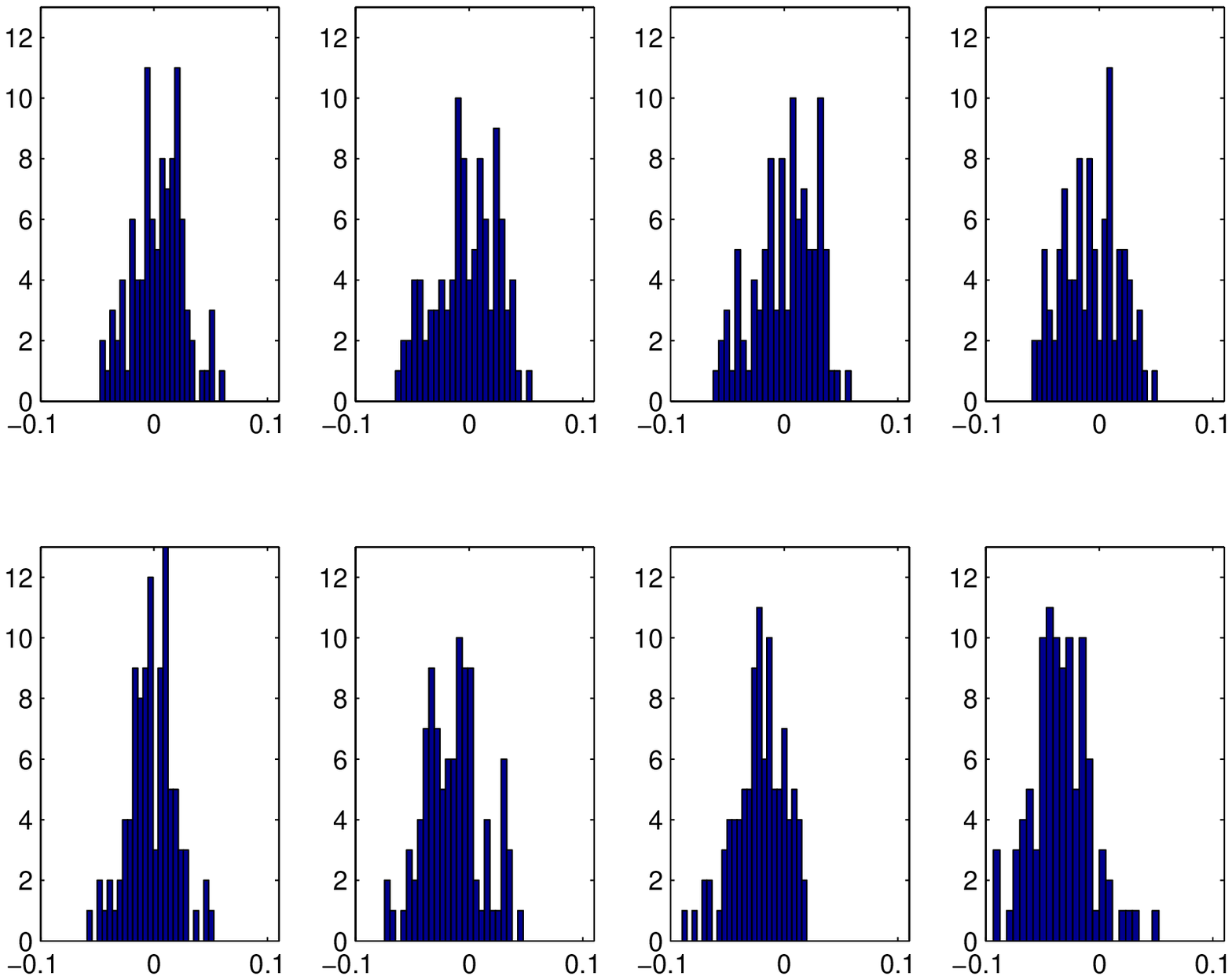}
\caption{Histograms for the estimation errors of  Efron's estimator for $\uo$ (top row) and $\hat{\mu}_0^*$ (bottom
  row). From left to
  right:  $\eps =  0.05$, $0.10$, $0.15$, and $0.20$. } 
\label{fig:uCompare}
\end{figure}

We now compare $(\hat{\sigma}_0^*,  \hat{\mu}_0^*)$  with the  estimators in
Efron \cite{Efron}.   Recall that one major difference between the two
approaches is that   Efron's estimators are not consistent for the
non-sparse case,    while  ours    are.   It is thus of interest to
make comparisons at different levels of sparsity.  To do so,    we set
$a$ at $1$, and   let  $\eps$ take four different values, $0.05, 0.10,
0.15$, and  $0.20$,  to represent different levels of
sparsity. For each $\eps$, we first generate samples according to   the main step  in the aforementioned experiment,    then   implement  $(\hat{\sigma}_0^*, \hat{\mu}_0^*)$   and the estimators of Efron \cite{Efron}, and finally repeat the experiment   for $100$ independent cycles.     The  results are   reported  in Figures \ref{fig:sCompare} - \ref{fig:uCompare}.    

The results show that our estimator of $\so^2$ is more
accurate than that of Efron \cite{Efron}, and the difference becomes
more prominent as $\eps$ increases.  
In fact,  when $\eps$ ranges between $0.05$ and $0.2$,   the estimation errors of $\hat{\sigma}_0^*$ are of the order $10^{-2}$, while those of Efron's estimator could get as large as the order $10^{-1}$.   
On the other hand,  the two estimators of $\mu_0$ are almost equally accurate, and the estimation errors for both approaches fluctuate around $0.02$ across different choices of $\eps$.   

 However,  the above comparison is  only for
 moderately large $n$.  With a much larger $n$,     the previous
 theory (Theorem \ref{thm:Main2b}) predicts  that  the estimation
 errors of  $(\hat{\sigma}_0^*, \hat{\mu}_0^*)$  will  become
 substantially smaller as  $(\hat{\sigma}_0^*, \hat{\mu}_0^*)$ is
 consistent for $(\so, \uo)$. In comparison,  the errors of Efron's estimators will not 
become  substantially smaller as  the estimators are not consistent.     To illustrate this point, we carry out a small scale
 simulation experiment. We take   $\eps = 0.1$ and $a = 1$ as 
 before,   while we let  $n  = 10^4$,  $4 \times 10^4$,  $1.6 \times 10^5$, and $6.4 \times 10^5$.
 For each $n$,   we generate samples according to  the main step, calculate the mean squared errors (MSE),   and repeat the process for $30$ independent cycles.  The results are
 reported in Table \ref{table:largen}, and they  support  the asymptotic analysis.   
\begin{table}[htb]
$$
\begin{tabular}{|c|l|c|c|c|c|}
\hline 
$n$ & & $10^4$ & $4 \times 10^4$ &$1.6 \times 10^5$ & $6.4 \times 10^5$\\
\hline
MSE for $\so$ & Efron's approach & 9.100  & 8.564 & 8.415 & 8.567  \\ 
\cline{2-6}
  &Our approach   & 0.816  & 0.276  & 0.047  & 0.031   \\
\hline\hline
MSE for $\uo$  &Efron's approach  & 8.916 &5.905   & 3.957  &3.617 \\  
\cline{2-6} 
&Our approach  &5.807  &3.019 & $1.106 $  &  0.538 \\
\hline
\end{tabular}
$$
\caption{Mean squared errors (MSE) for various values of $n$.  
The corresponding  MSE equals the value in each cell  times $10^{-4}$.  }
\label{table:largen}
\end{table} 

\begin{figure}[htb]
\centering
\includegraphics[height = 1.75 in, width = 4  in]{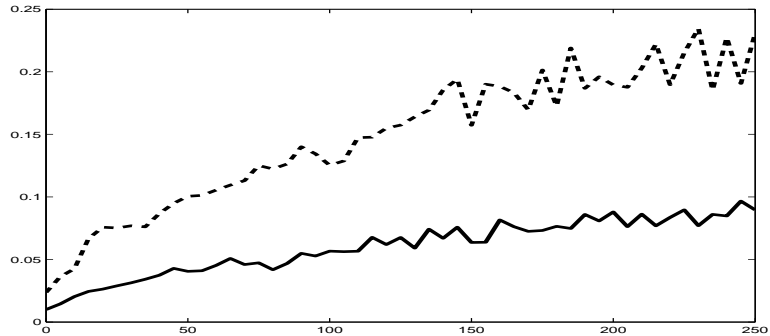}
\caption{$x$-axis:  $L$.  $y$-axis:   root mean squared error for $\hat{\mu}_0^*$ (dashed) and $\hat{\sigma}_0^*$ (solid).}
\label{fig:depen}
\end{figure}

Finally, we investigate the performance of the proposed procedures for
dependent data. Fix $n  = 10^4$, $\eps = 0.1$, and $a = 1$, and let $L$
range from $0$ to $250$ with an increment of $5$. For each $L$,
generate $n + L$ samples $w_1, w_2, \ldots, w_{n+L}$ from $N(0,1)$
  and  let  $z_j =  (\sum_{k = j}^{k= j + L} w_k)/\sqrt{L + 1}$,
so that $\{z_j\}_{j=1}^n$  are  block-wise dependent (block size  equal
to  $L+1$) and  the marginal distribution of each $z_j$ is  $N(0,1)$. At the same time, generate the mean vector $\mu$ and the
vector of standard deviations $\sigma$ according to the main step,  
let $X_j = \mu_j + \sigma_j  \cdot  z_j$,     and implement $(\hat{\mu}_0^*,
\hat{\sigma}_0^*)$ to $\{X_j\}_{j=1}^n$. We then repeat the process
for $100$ independent cycles. The results are reported in Figure
\ref{fig:depen}, which suggests that the estimation errors increase as
the range of dependency increases. However, when $L \leq 100$, for
example, the  estimation errors  are still relatively small,
especially those for $\sigma_0^*$. This suggests that the procedures
are relatively robust to short range dependency.   

\section{Applications to microarray analysis}  
\label{sec:Appli}
\setcounter{equation}{0}

We now apply the proposed  procedures to the analysis of 
  the breast cancer
and HIV microarray data sets that were analyzed in Efron \cite{Efron}. 
The R code for our procedures is available on the web at 
{\it http://www.stat.purdue.edu/$\tilde{\;}$jinj/Research/software}.
The $z$-scores for  both   data sets   
 can   be downloaded from this site as well; they were kindly  provided by Bradley Efron.
The R code for  Efron's procedures and related software can be
downloaded from  {\it  http://cran.r\-project.org/src/contrib/Descriptions/ locfdr.html}.  
For reasons of space, we focus on the  breast cancer data   
and only comment briefly on the  HIV data.  

The breast cancer data was based on $15$  patients diagnosed with
breast cancer,  $7$ with the BRCA1 mutation and $8$ with the BRCA2
mutation. Each patient's tumor was analyzed on a separate microarray,   
and the microarrays reported on the same set of $N = 3226$ genes.  For  the $j$-th gene,  
the two-sample $t$-test comparing the seven BRCA1 responses  with the eight 
BRCA2 was computed.    The $t$-score  $y_j$  was  first  converted to the 
$p$-value by   $p_j   =  \bar F_{13}(y_j)$,   and   was then  converted to the $z$-scale \cite{Efron}, $X_j  =  \bar \Phi^{-1}(p_j)  = \bar \Phi^{-1}(\barF_{13}(y_j))$,
where $\bar\Phi$ and $\bar F_{13}$ are the survival functions of
$N(0,1)$  and $t$-distribution with $13$ degrees of
freedom,  respectively. 

We model $X_j$ as $N(\mu_j, \s_j^2)$ variables
with weakly   dependent structure, and for a pair of unknown
parameters $(\uo, \so)$, $(\mu_j, \sigma_j) = (\uo, \so)$ if and only
if the $j$-th gene is not differentially expressed. Since $X_j$ is
transformed from the $t$-score which has been standardized by
the corresponding standard error, it is reasonable to assume that the  
null effects are homogeneous, and that  all effects are homoscedastic;
see for example, \cite{Churchill, Efron}.   
The normality assumption is also reasonable here, as the
marginal density of non-null effects can generally be well
approximated by Gaussian mixtures; see  \cite[Page
99]{Efron}.  Particularly, it is well known that  the set of all
Gaussian mixing densities is dense in the set of all density
functions under the $\ell^1$-metric.   
 
We now proceed with  the data analysis.  
The analysis includes  three parts:  estimating the null parameters
$(\so, \uo)$, estimating the proportion of  non-null effects, and
implementing the local FDR  approach proposed by
Efron et al. \cite{Efronetal}.  

The first part is estimating  $(\so, \mu_0)$.   We apply
$(\hat{\sigma}_0^*, \hat{\mu}_0^*)$ (defined in
(\ref{EqDefinemusigmahat*}))  as well as the estimators used by Efron \cite{Efron} to the $z$-scores.   For the breast cancer data,   our procedure 
yields   $(\hat{\sigma}_0^*, \hat{\mu}_0^*) = (1.5277,    -0.0525)$,  
while Efron's estimators give 
$(\hat{\sigma}_0, \hat{\mu}_0)  = (1.616,   -0.082)$.

The second part of the analysis is estimating the proportion of  
non-null effects.  We implement our procedure as well as Meinshausen
and Rice's \cite{Rice} approach and the approach of Cai et al.  
\cite{CJL} (which we denote by MR and CJL respectively for short), to the  $z$-scores of the breast cancer data.    The bounding function
$a_n^*$ for MR estimator is set as $1.25
\times \sqrt{2 \log \log n}/\sqrt{n}$,  and the $a_n$ for CJL  estimator   is   set as $\sqrt{2 \log \log n}/\sqrt{n}$;   see \cite{CJL} for details.     Using  the estimated null parameters either obtained  by  Efron's approach or obtained by our  approach,    we apply each of these  procedures   to the $z$-scores.       In addition, the local FDR  approach  also provides an estimate for the proportion automatically. The results   are reported in Table
\ref{table:esteps}.

\begin{table} 
$$
\begin{tabular}{|c|c|c|c|c|}
\hline 
   &Our estimator &Local FDR & MR &CJL \\
\hline
Our Estimated Null  & 0.0040  &0.0128    &  0.0033  &  0    \\ 
\hline
 Efron's  Estimated Null    & 0  & 0  &  0.0098    & 0    \\
\hline
\end{tabular}
$$
\caption{Estimated proportion of non-null effects for the breast
  cancer data.}
\label{table:esteps}
\end{table}

In the last part of the analysis we implement the local FDR
thresholding procedure proposed in \cite{Efronetal} with the $z$-scores of the breast cancer data.   For any given FDR-control parameter $q \in (0,1)$,  
the procedure calculates a score for each data point and  determines a
threshold $t_q$ at the same time. A hypothesis is rejected
if the score exceeds the threshold and is accepted otherwise.   
If we call a rejected hypothesis   a ``discovery,''  then the local FDR thresholding procedure  controls  the expected false discovery rate at level $q$, 
$E [\frac{\# \mbox{False Discoveries}}
{\# \mbox{Total Discoveries}}]  \leq q$.   
See \cite{Efronetal} for details.  

With Efron's
estimated  null parameters,  for any  fixed $q \in (0,1)$,    the  local FDR
procedures report {\it no} rejections for the breast cancer data set. Also,    three different estimators 
 for the proportion report  $0$.  
These suggest  that either the proportion of signals (differentially
expressed genes) is small and/or the signal is very weak.   

In contrast, with our estimated null parameters, the estimated proportions  are small but  nonzero.    Furthermore,   the local FDR procedures report rejections when  $q \geq 0.91$.  For example,   the number of total discoveries equal to $167$ when   $q = 0.92$, and equal to $496$ when  $q = 0.94$.  Take $q = 0.94$,    for example,    since for any   $q \in (0,1)$, the number of true discoveries approximately  equal to $(1-q)$ times the number of total discoveries \cite{Efron},   this suggests  a total of $30$ true discoveries.    
The  result  is  consistent with biological discoveries.    Among the $496$ genes which are identified 
to be differentially expressed by the local FDR procedures, $17$ of
them have been discovered in the study by Hedenfalk et al.
\cite{BreastCancer}. The corresponding  Unigene cluster IDs are:
 Hs.182278, Hs.82916, Hs.179661, Hs.119222,
Hs.10247, Hs.469, Hs.78996, Hs.11951, Hs.79078, Hs.9908,
Hs.5085, Hs.171271, Hs.79070,   Hs.78934, Hs.469,
Hs.197345, Hs.73798.
We also identified several genes whose
functions are associated with the cell cycle, including PCNA, CCNA2,
and CKS2. These genes are found to be significant by Storey  et al.  
\cite{Storey2}. The results  indicate that our estimated null
parameters lead to reliable  identification of   differentially expressed
genes.

Similarly,   for  the   HIV data,  our estimators give  
$(\hat{\sigma}_0^*, \hat{\mu}_0^*) = (0.7709, -0.0806)$,  
while Efron's method gives 
$(\hat{\sigma}_0, \hat{\mu}_0)     =   (0.738,   -0.082)$.  With $q = 0.05$,  the local FDR procedures report $59$ total discoveries with our estimated null parameters,   and $80$ with Efron's estimated null parameters;  the latter yields slightly more signals.   

\section{Proofs of the main results} 
\label{sec:Proof}
\setcounter{equation}{0}

We now prove  Theorems \ref{thm:Main1},
\ref{thm:Main2b},  and \ref{thm:esteps}.  The proof of Theorem
     \ref{thm:estepsAdd} is similar to those of  Theorems  \ref{thm:Main2b} and
\ref{thm:esteps} and so is omitted. As the proofs for the estimators
of $\so^2$  and $\uo$ are similar, we focus on $\so^2$.
We first collect a few technical results and outline the basic ideas. 
The proofs of these preparatory lemmas are given in \cite{JC}.
      

\begin{lemma}  \label{lemma:Est1}
Let $\so^2(\cdot; \cdot)$ and $\uo(\cdot; \cdot)$ be defined as in
(\ref{EqDefineso}).   Fix $t > 0$. For any  differentiable complex-valued  functions $f$ and $g$ satisfying $|f(t)| \neq 0$ and $|g(t)| \neq 0$,  
\[
|\so^2(f,t) - \so^2(g,t)| \leq \frac{|g(t)|}{t |f(t)|^2}    \bigl[ \bigl(2 t \cdot |\so^2(g,t)|  +  |\frac{g'(t)}{g(t)}| \bigr)   |f(t)  - g(t)|  +   |f'(t) - g'(t)| + r_n^{(1)}(t)  \bigr],   
\] 
\[
|\uo(f,t) - \uo(g,t)| \leq   \frac{|g(t)|}{|f(t)|^2} \cdot  \bigl[\bigl(2|\uo(g,t)|  + |\frac{g'(t)}{g(t)}| \bigr) \cdot |f(t)  - g(t)|   + |f'(t) -g'(t)|  + r_n^{(2)}(t) \bigr],  
\]
where 
$
r_n^{(1)}(t) =  \frac{1}{|g(t)|} \cdot \bigl[t \cdot |\so^2(g,t)| \cdot   |f(t) - g(t)|^2 + |f(t) - g(t)| \cdot |f'(t) -g'(t)| \bigr]   
$
 and  
$
 r_n^{(2)}(t) = \frac{1}{|g(t)|}\cdot \bigl[ |\uo(g,t)| \cdot  |f(t) - g(t)|^2 + |f(t)  - g(t)| \cdot |f'(t) - g'(t)| \bigr].
$
\end{lemma}
Heuristically,   $|\varphi(\that_n)|/|\varphi_n(\that_n)|^2  \sim n^{\gamma}$,   $\so^2(\varphi, \that_n)  \sim   \so^2$,   $|\varphi'(\that_n)|/|\varphi(\that_n)|   \sim  \so^2 \that_n$,  and  
\begin{equation}  \label{EqAdd101}
|\varphi_n(\that_n) - \varphi(\that_n)|   \leq O_p(\sqrt{\log n}/\sqrt{n}), \qquad   |\varphi_n'(\that_n) - \varphi'(\that_n)|   \leq O_p(\sqrt{\log n}/\sqrt{n}).
\end{equation} 
Applying  Lemma \ref{lemma:Est1} with $f = \varphi_n$, $g = \varphi$,  and $t  =  \that_n(\gamma)$, we have   \begin{align*}
&\qquad  |\so^2(\varphi_n, \that_n(\gamma)) - \so^2(\varphi, \that_n(\gamma))|   \\  
&\sim    n^{\gamma}  \bigl(3    \so^2       |\varphi_n(\that_n(\gamma)) - \varphi(\that_n(\gamma))|  + \frac{1}{\that_n(\gamma)} |\varphi_n'(\that_n(\gamma)) - \varphi'(\that_n(\gamma))|\bigr)  
\sim O(n^{\gamma -\frac{1}{2}} \sqrt{\log n}),  
\end{align*}
and Theorem \ref{thm:Main1} follows.   We now study  (\ref{EqAdd101}) in detail.   
\begin{lemma} \label{lemma:phi}
Set
$
W_0(\varphi_n; n)   =   W_0(\varphi_n; n, X_1, \ldots, X_n)  =  \sup_{0 \leq t \leq \log n}  |\varphi_n(t)  - \varphi(t)|.
$
Fix $q_1 > 3$. Let $\Lambda_n(q,A; \uo, \so, \eps_0)$  be given as  in  Theorem \ref{thm:Main1}.   When $n \goto \infty$,   
\[
\sup_{\{(\mu, \sigma) \in \Lambda_n(q,A; \uo, \so, \eps_0)\}} P\{W_0(\varphi_n; n) \geq  \sqrt{2 q_1  \log n}/\sqrt{n}  \}  \leq  4 \log^2(n) \cdot   n^{-  q_1/3}  \cdot (1 + o(1)).  
\]
\end{lemma}

Lemma \ref{lemma:phi} implies that except for  an event with
algebraically small probability,   
$|\varphi(\that_n)  - \varphi(t_n)|  \leq  W_0(\varphi_n; n) \leq \sqrt{2 q_1 \log n} /\sqrt{n}$.   
This naturally leads to a precise description of the stochastic  behavior of 
  $|\that_n(\gamma) - t_n(\gamma)|$ given in the following lemma. 
\begin{lemma}  \label{lemma:t}
Let  $q_1 > 0$ and let $\Lambda_n(q,A; \uo, \so, \eps_0)$,
$\that_n(\gamma)$,  and $t_n(\gamma)$  be given as   in  Theorem \ref{thm:Main1}.    When $n \goto \infty$, 
\[
\sup_{\{(\mu, \sigma) \in \Lambda_n(q,A; \uo, \so, \eps_0)\}}   \{ |\that_n(\gamma)  - t_n(\gamma)|  \cdot 1_{\{W_0(\varphi_n; n) \leq \sqrt{2 q_1 \log n}/\sqrt{n}\}} \}  \leq  \frac{1}{\so}  \sqrt{\frac{q_1}{\gamma}}    n^{\gamma - 1/2}   (1 + o(1)).  
\]
\end{lemma} 

We now  study $|\varphi_n'(\that_n) - \varphi'(\that_n)|$.   Pick a
constant  $\pi_0 >  \frac{1}{\so}   \sqrt{q_1/\gamma}$ and set 
\[
W_1(\varphi_n, \gamma,  \pi_0; n)   =   W_1(\varphi_n, \gamma, \pi_0; n, X_1, \ldots, X_n)  =  \sup_{|t -  t_n(\gamma)| \leq  \pi_0 \cdot  n^{\gamma - 1/2}}  |\varphi'_n(t)  - \varphi'(t)|.    
\]
By Lemma \ref{lemma:t},   except for an event with  algebraically small  probability,  
$|\that_n(\gamma)  - t_n(\gamma)|  \leq     \pi_0  \cdot   n^{\gamma - 1/2}$,   
and consequently  
$|\varphi_n'(\that_n(\gamma)) - \varphi'(\that_n(\gamma))|  \leq  W_1(\varphi_n, \gamma,  \pi_0; n)$.    
The following lemma describes the tail behavior of $W_1$. 
\begin{lemma} \label{lemma:phi'} 
Fix  $\gamma \in (0,1/2)$, $\pi_0 >  \frac{1}{\so}
\sqrt{q_1/\gamma}$ and set $\bar{s}_n^2  =  \frac{1}{n} \sum_{j = 1}^n   E[X_j^2]$.    There exist  constants $C_1$ and  $C_2 > 0$ such that  for any $(\mu, \sigma) \in \Lambda_n(q,A; \uo, \so, \eps_0)$,  
$\bar{s}_n  \leq C_1$,  
\[
P\{W_1(\varphi_n,  \gamma, \pi_0; n)  \geq  \bar{s}_n  \cdot \frac{\sqrt{(q -2) \log n} + 2 \bar{s}_n}{\sqrt{n}}  \}   \leq    C_2 \cdot n^{-c_1(q,\gamma)},   
\]
where $c_1(q,\gamma)$ is as in Theorem \ref{thm:Main1}. 
As a result,  except for an event with algebraically small probability,  
$|\varphi'_n(\that_n(\gamma))  - \varphi'(\that_n(\gamma))| \leq  W_1(\varphi_n, \gamma, \pi; n)  \leq   O(\sqrt{\log n}/\sqrt{n})$. 
\end{lemma} 

We have now elaborated the inequalities in (\ref{EqAdd101}).   The only missing piece is the following lemma,  
which  gives the basic properties of $\so^2(\varphi; t)$ and $\uo(\varphi; t)$. \begin{lemma}  \label{lemma:g}
Fix   $q \geq 3$ and  $A > 0$,     with   $\psi(t)$ and $\tau_n$ as defined in  (\ref{EqDefinepsi}) and (\ref{EqDefinetaun}) respectively,  write $\psi(t) = \eps_n g(t)$ and $r(t) = \frac{\eps_n}{1-\eps_n} r(t)$. For all   $(\uo, \sigma_0, \eps_0)$-eligible $(\mu, \sigma)$ and all $t > 0$,      there is a constant $C > 0$ such that
\begin{align}
|\so^2(\varphi,t) - \so^2|   &\leq    \frac{|r'(t)|}{t}  \cdot \frac{1 +   |r(t)|}{|1 +   r(t)|^2} \leq C |\psi'(t)|/t,  \label{lemmag1} \\ 
|\uo(\varphi,t) - \uo|   &\leq     |r'(t)|  \cdot \frac{1 +     |r(t)|}{|1 + r(t)|^2}  \leq C|\psi'(t)|.     \label{lemmag2} 
\end{align}
Additionally,       uniformly for  all  $(\mu, \sigma)\in
\Lambda_n(q,A; \uo, \so, \eps_0)$ and all $t > 0$,    
\begin{itemize}
\item[(a1).] $|g(t)| \leq e^{-\frac{\tau_n^2
      t^2}{2}} \leq 1$,   $|g'(t)|  \leq A$,   $|g''(t)| \leq C(1 +
  A^2)$,  $|g'''(t)|\leq C(1 + A^3)$, and  $|g'(t)|  \leq  A e^{-
    \frac{\tau_n t^2}{2}}  +  \mmin \{  A^2 t e^{-\frac{\tau_n
      t^2}{2}}, \frac{2}{et} \}$;   
\item[(a2).]   consequently,        $|\varphi'(t)|/|\varphi(t)|  =
  \so^2 \cdot t  \cdot (1  + o(1))$;  
\item[(a3).]   the second derivative of   $\so^2(\varphi; t)$ is
  uniformly bounded, and $\so^2(\varphi; t) \goto
  \so^2$, $\frac{d}{dt} \so^2(\varphi; t) \goto 0$ as $t \goto
  \infty$.     
\end{itemize}
Similarly,   both $\uo(\varphi;t)$ and  its first two derivatives are uniformly bounded for all $t > 0$, and $\frac{d}{dt}\uo(\varphi;t) \goto 0$ if $\uo(\varphi;t) \goto \uo$.     
\end{lemma} 
 
We now prove Theorem \ref{thm:Main1}, \ref{thm:Main2b}, and \ref{thm:esteps}. 

\noindent
{\bf Proof of Theorem \ref{thm:Main1}:}
  Since the arguments are  similar, we  prove
the first claim only.    Write $\that_n  = \that_n(\gamma)$,    $t_n =
t_n(\gamma)$,  and  $W_1(\varphi_n; n) = W_1(\varphi_n, \gamma, \pi_0;
n)$.  Pick constants $q_1$ and $\pi_0$ such that   $1 <  q_1/
\mmax\{3,  (q - 1 - 2 \gamma) \}   <     2$ and $\pi_0 >
\frac{1}{\sigma_0} \sqrt{q_1/\gamma}$. Introduce  events
\[
B_0    =   \{W_0(\varphi_n;n)   \leq   \sqrt{2 q_1 \log n}\}, \;\;     B_1   = \{W_1(\varphi_n;  n)   \leq  \frac{s_n   \sqrt{(q -2) \log n} + 2 s_n^2}{\sqrt{n}}\}. 
\]
Note that the choice of $q_1$ satisfies
$c_1(q,\gamma) < q_1/3$ and $c_2(\so, q, \gamma)  >  \so^2 \sqrt{2 q_1}$,  
where  $c_1(q,\gamma)$ and $c_2(\so,q,\gamma)$ are defined as  in Theorem \ref{thm:Main1}. 
Use Lemma \ref{lemma:phi}  and Lemma \ref{lemma:phi'},     $P\{B_0^c\}  \leq   \bar{o}(n^{-  q_1/3})$ and 
$P\{B_1^c\}  \leq  \bar{o}(n^{- c_1(q, \gamma)})$;  since $c_1(q, \gamma) < q_1/3$,       $P\{B_0^c  \cup  B_1^c \}   \leq    \bar{o}(n^{-c_1(q,\gamma)})$.    We now focus on $B_0 \cap B_1$.   By triangle inequality,   
$|\so^2(\varphi_n; \that_n) - \so^2(\varphi; t_n)|  \leq
|\so^2(\varphi_n; \that_n) - \so^2(\varphi; \that_n)| +
|\so^2(\varphi; \that_n) - \so^2(\varphi; t_n)|$. 
Note that by the choice of $\pi_0$ and  Lemma \ref{lemma:t},  
$|\that_n   -   t_n|  \leq  \pi_0  \cdot n^{\gamma - 1/2}$ 
for sufficiently large $n$,  it thus  follows from
Lemma \ref{lemma:g} that $|\so^2(\varphi; \that_n) -
\so^2(\varphi; t_n)|\sim o(|\that_n - t_n|) = o(n^{\gamma -1/2})$;  recall    $c_2(\sigma_0, q,\gamma) > \so^2\sqrt{2 q_1}$,  so to show the claim,   it suffices to show that as $n \goto \infty$,   
\begin{equation}   \label{EqthmMain1ToShow1}
|\so^2(\varphi_n; \that_n)  - \so^2(\varphi; \that_n)|  \leq  3  \so^2  \cdot  \sqrt{2 q_1 \log n}  \cdot n^{\gamma - 1/2}  \cdot (1 + o(1)),   \qquad \mbox{over  $B_0 \cap B_1$}.   
\end{equation} 

We now show (\ref{EqthmMain1ToShow1}).   Over the event $B_0 \cap
B_1$,  recall $|\that_n - t_n| \leq \pi_0 n^{\gamma-1/2}$,   so by  (\ref{EqDefinetnadd}),     $\that_n  \sim t_ n   \sim  \sqrt{2 \gamma \log n}/\so$;  by  Lemma \ref{lemma:g},  this implies  $\so^2(\varphi, \that_n) \sim \so^2$ and  $|\varphi'(\that_n)|/|\varphi(\that_n)|   \sim  \so^2 \that_n  \sim \so^2 t_n$.  
Moreover,  since  $|\varphi_n(\that_n)  - \varphi(\that_n)| \leq \sqrt{2 q_1 \log  n}/\sqrt{n}$,  it follows that   
$|\varphi(\that_n)|/(\that_n |\varphi_n(\that_n)|^2)  \sim  (1/t_n) n^{\gamma}$.   Lastly,   by Lemma \ref{lemma:phi'},       $|\varphi_n'(\that_n)  - \varphi'(\that_n)| \leq  O(\sqrt{\log n}/\sqrt{n})$.       Combining these,   (\ref{EqthmMain1ToShow1}) follows directly by applying 
  Lemma \ref{lemma:Est1} with   $f = \varphi_n$, $g =
\varphi$, and $t = \that_n$.    
  \qed


\medskip\noindent
{\bf Proof of Theorem \ref{thm:Main2b}:}
Note that,   by triangle inequality, $|\so^2(\varphi_n;
\that_n(\gamma)) - \so^2| \leq |\so^2(\varphi_n; \that_n(\gamma)  -
\so^2(\varphi; t_n(\gamma))| + |\so^2(\varphi; t_n(\gamma))  -
\so^2|$. Theorem \ref{thm:Main2b} now follows directly from Theorem
\ref{thm:Main1} and Lemma \ref{lemma:g}. \qed

\medskip\noindent
{\bf Proof of Theorem \ref{thm:esteps}:}
Without loss of generality, set $\uo = 0$ and $\so = 1$.  
Write   $t_n  = \sqrt{2 \gamma \log n}$,  $\eps_n =
\eps_n(\mu, \sigma)$,    $\varphi_n(t) = \varphi_n(t; X_1, \ldots,
X_n, n)$,   $\varphi(t) = \varphi(t;  \mu, \sigma, n)$,
$\Omega_n(t) =  \Omega_n(t; X_1, \ldots, X_n, n)$,  and $\Theta_n =
\Theta_n(\gamma; q,A,\uo, \so,\eps_0)$.    Set $\Omega(t) =
E[\Omega_n(t)]$,
$\Psi_n^*(t) =   \sup_{\{0 \leq s  \leq t \}}  \{1 -  \Omega_n(s)\}$,
and  
$\Psi^*(t) =   \sup_{\{0 \leq s  \leq t \}}  \{1 - \Omega(s)\}$. 
Note that it is sufficient to show that when $n \goto \infty$, 
 (a)    except for an event with algebraically  small probability, 
$\sup_{\{ (\mu, \sigma) \in \Theta_n \}}  |\Psi_n^*(t_n) -     \Psi^*(t_n)|  \leq O( \log^{-3/2}(n)  \cdot   n^{\gamma -1/2})$,    and  (b)     $\sup_{\{(\mu, \sigma) \in  \Theta_n \}}  |\frac{\Psi^*(t_n)}{\eps_n}  - 1|  =  o(1)$.  

We first show (a). By symmetry, $|\Psi_n^*(t_n) - \Psi^*(t_n)|$ does not exceed
\begin{equation}  
   \sup_{0 \leq t \leq t_n}   |\Omega_n(t) - \Omega(t)|     
\leq 2  \int_{0}^1 (1 - \xi)  e^{\frac{t_n^2 \xi^2}{2}} \sup_{0 \leq t \leq t_n}  |\mathrm{Re}(\varphi_n(t)) - \mathrm{Re}(\varphi(t))|   d\xi.   \label{thm:estepsPf11}
\end{equation} 
Moreover, similar to the proof of Lemma 7.2 in \cite{Jin},   we have that for fixed $q > 3/2$,  
$\sup_{\{(\mu, \sigma) \in \Theta_n \}} \sup_{\{0 \leq t \leq t_n\}} |\mathrm{Re}(\varphi_n(t)) - \mathrm{Re}(\varphi(t))|   \leq O(\sqrt{\log n}/\sqrt{n})$  
except for an event with probability $\sim  2 \log^2(n) \cdot n^{-2q/3}$.
Elementary calculus  yields
$|\Psi_n^*(t_n) - \Psi^*(t_n)|    \leq  O(\sqrt{\log n}/\sqrt{n}) \cdot  \int_0^1(1 - \xi) e^{(\gamma \log n) \cdot   \xi^2} d \xi   
 = O(\log^{-3/2}(n)   \cdot  n^{\gamma - 1/2})$,
and (a) follows.

We now show (b). Let  $\hat{f}$ be the Fourier transform of $f$ and
let $\phi_{\delta_j(t)}(x)$ be the density function of $N(0,
\delta_j^2(t))$ with  $\delta_j(t)   =   t (\sigma_j^2-1)^{1/2}$. Set
$\rho(x) = 2(1 - \cos(x))/x^2 $ for $x \neq 0$ and $\rho(0) = 1$.
Elementary calculus shows that
$
\hat{\phi}_{\delta_j(t)}(\xi)  = \mathrm{exp}(\frac{(1 -
  \sigma_j^2)t^2 \xi^2}{2})$ and  $\hat{\rho}(\xi) =  \mmax\{1 -
|\xi|,  0\}.
$
So by the Fourier Inversion Theorem \cite[Page 22]{Mallat}, 
\begin{align*}
\Omega(t)   &=\frac{1}{n}  \sum_{j = 1}^n  \int_{-1}^1 (1 - |\xi|)   \mathrm{exp}(\frac{(1 - \sigma_j^2) t^2 \xi^2}{2})  \cos(t    \mu_j \xi ) d\xi   \\
&=   \frac{1}{n}  \sum_{j = 1}^n  \int_{-1}^1  \hat{\phi}_{\delta_j(t)}(\xi)  \hat{\rho} (\xi)  \cos(t  \mu_j \xi ) d \xi  =   \frac{1}{n}  \sum_{j = 1}^n  \phi_{\delta_j(t)}*\rho(t \mu_j ), 
\end{align*}
 where $*$ is the usual convolution. Since $\phi_{\delta_j(t)}*\rho(t \mu_j) = 1$ when $(\mu_j, \sigma_j) = (0,1)$,    
\begin{equation}  \label{thm:estepsPf4}
1 - \Omega(t)  =    \eps_n   \cdot   \mathrm{Ave}_{\{j:\; (\mu_j,\sigma_j) \neq (0,1)\}}  \{ 1 -  \phi_{\delta_j(t)}*\rho(t \mu_j )\}.  
\end{equation} 
Note that  $\phi_{a_n}*\rho(b_n) \goto 0$    for any sequences $\{a_n\}_{n =1}^{\infty}$ and
$\{b_n\}_{n=1}^{\infty}$  satisfying  $\mmax\{a_n, b_n\} \goto \infty$,  so by (\ref{thm:estepsPf4}) and the definition of $\Theta_n$,
$ \sup_{\{(\mu, \sigma) \in  \Theta_n \}}  \big|\frac{1  - \Omega(t_n)}{\eps_n} - 1 \big|  =  o(1)$.  Note that $0 \leq \phi_{\delta_j(t)}*\psi(t) \leq 1$ for all $t$,   so by (\ref{thm:estepsPf4}) and the definition of $\Psi^*$,      $\Omega(t_n) \leq    \Psi^*(t_n)      \leq \eps_n$; as a result,     $\big|\frac{1  - \Psi^*(t_n)}{\eps_n} - 1 \big|   \leq \big|\frac{1  - \Omega(t_n)}{\eps_n} - 1 \big|$,    and (b) follows directly.   \qed

\subsection*{Acknowledgments}
We thank Bradley Efron  for references and
kindly sharing the data sets. We thank Paul Shaman for a careful reading of our manuscript and for suggestions which lead to significant improvement of the presentation of the paper.   We also thank Herman Rubin, an Associate editor, and referees for helpful comments and references.    


\section{Appendix}  
\label{sec:appendix}
\setcounter{equation}{0}
\subsection{Proof of Theorem \ref{thm:Depen}}
For short, write $\that_n  = \that_n(\gamma)$ and $t_n = t_n(\gamma)$.  
The following  two lemmas are proved  in Section \ref{subsubsec:lemmaDepen1} and Section  \ref{subsubsec:lemmaDepen2}  respectively.    
\begin{lemma}  \label{lemma:Depen1}
With $\alpha(\cdot)$ and  $\tilde{\Lambda}_n(a,B,q,A)$ as in Theorem \ref{thm:Depen}.   Fix  $r \in (1.5,  d - (2d + 2.5)\gamma)$.       As $n \goto \infty$,  uniformly for all $(\mu, \sigma) \in \tilde{\Lambda}_n(a,B,q,A)$,    except for an event  with a probability of $\bar{o}(n^{-2r/3})$,  
$\sup_{\{0 \leq t \leq \log n\}}  |\varphi_n(t)  - \varphi(t)|   =   \bar{o}(n^{-(d-r)/(2d + 2.5)})$.  
\end{lemma} 
\begin{lemma}  \label{lemma:Depen2}
With $\alpha(\cdot)$ and  $\tilde{\Lambda}_n(a,B,q,A)$ as in Theorem \ref{thm:Depen}. Fix $\gamma \in (0, \frac{d -1.5}{2d + 2.5})$ and an integer $k \geq 0$.     As $n \goto \infty$,   for all $(\mu, \sigma) \in \tilde{\Lambda}_n(a,B,q,A)$,   $\sup_{\{0 \leq t \leq \log n\}} \{|\varphi_n^{(k)}(t)  - \varphi^{(k)}(t)| \}   \leq   O_p(\log^{ak}(n))$, and 
 $[\varphi_n^{(k)}(t_n)  - \varphi^{(k)}(t_n)]= O_p(\log^{(a + 1/2)k}(n)/\sqrt{n})$.  
\end{lemma}

To show the theorem,   it is sufficient to show that 
\begin{equation} \label{EqlemmaDepen1}
|\varphi_n(\that_n) - \varphi(\that_n)|  = O_p(1/\sqrt{n}), \qquad   |\varphi'_n(\that_n) - \varphi'(\that_n)|  = O_p(\log^{a + 1/2}(n)/\sqrt{n}).   
\end{equation} 
In fact, by  triangle inequality,  
\begin{equation}  \label{EqlemmaDepen1Add} 
|\so^2(\varphi_n; \that_n)  - \so^2(\varphi; t_n)|   \leq  |\so^2(\varphi_n; \that_n)  - \so^2(\varphi; \that_n)|    +  |\so^2(\varphi; \that_n)  - \so^2(\varphi; t_n)|.      
\end{equation} 
Once (\ref{EqlemmaDepen1}) is proved,        by   similar arguments  as in the proof of  Theorem \ref{lemma:t}, 
\begin{equation}  \label{EqlemmaDepenAdd101} 
|\that_n   -  t_n|  =   O_p(n^{\gamma -1/2}),  
\end{equation} 
it thus follows from Lemma \ref{lemma:g} that 
\begin{equation}  \label{EqlemmaDepen4}
|\so^2(\varphi; \that_n)  - \so^2(\varphi; t_n)|  =  o_p(|\that_n - t_n|) =  o_p(n^{\gamma - 1/2}). 
\end{equation}  
At the same time,   
by (\ref{EqlemmaDepenAdd101})  and Lemma \ref{lemma:g},   except for an event with   asymptotically vanishing probability,  
$|\varphi_n(\that_n)|/|\varphi_n(\that_n)|^2 \sim n^{\gamma}$,  $\so^2(\varphi; \that_n) \sim \so^2$, and   $|\varphi'(\that_n)|/|\varphi(\that_n)| \sim \so^2 \that_n$;    
applying  Lemma \ref{lemma:Est1} with $f = \varphi_n$, $g = \varphi$, and $t = \that_n$,   it follows that
\begin{equation}  \label{EqlemmaDepen3}
|\so^2(\varphi_n; \that_n)  - \so^2(\varphi; \that_n)|  = O_p(n^{\gamma - 1/2}). 
\end{equation} 
   The theorem follows directly by  inserting   (\ref{EqlemmaDepen4})  and  (\ref{EqlemmaDepen3}) into     (\ref{EqlemmaDepen1Add}).   

We now show (\ref{EqlemmaDepen1}).  Since  the proofs are similar, we only show the first equality.   
Applying Lemma \ref{lemma:Depen1} with   $r =  (1.5 +  d - (2d + 2.5) \gamma)/2$, it follows that there is an event $A_n$ such that $P\{A_n^c\}$ is  algebraically small  and 
\begin{equation}  \label{EqlemmaDepen5} 
\sup_{\{0 \leq t \leq \log n\}}   |\varphi_n(t) -  \varphi(t)| \leq   \bar{o}(n^{-\frac{1}{2}\cdot (\gamma  + \frac{d - 1.5}{2d + 2.5})}), \qquad \mbox{over $A_n$}.          
\end{equation}   
By similar arguments as in the proof of Lemma \ref{lemma:t},   it follows that   
\begin{equation}  \label{EqlemmaDepen6}
|\that_n(\gamma)   -  t_n(\gamma)|    \leq  \bar{o}(n^{\frac{1}{2}\cdot   (\gamma  - \frac{d - 1.5}{2d + 2.5})}), \qquad \mbox{over $A_n$},  
\end{equation}  
notice the exponent is negative.   Now,  let   $\ell$ be the smallest integer  satisfying  $(\ell+1)  \cdot  |\gamma -   \frac{d - 1.5}{2d + 2.5}| > 1$.     By Taylor expansion,  for some $\xi$ falling between $\that_n$ and $t_n$,   
\[
\varphi_n(\that_n) - \varphi(t_n)    =  \sum_{k = 0}^{\ell}   \frac{\varphi_n^{(k)}(t_n) - \varphi^{(k)}(t_n)}{k!}  (\that_n - t_n)^{k}  +   \frac{\varphi_n^{(\ell+1)}(\xi) - \varphi^{(\ell+1)}(\xi)}{(\ell + 1)!}     (\that_n - t_n)^{\ell+1}. 
\]  
Notice that by the  choice of $\ell$ and (\ref{EqlemmaDepen6}),       $(\that_n - t_n)^{\ell+1}  =   \bar{o}(1/\sqrt{n})$ over $A_n$,  the claim follows directly from  Lemma \ref{lemma:Depen2}.     \qed

\subsubsection{Proof of Lemma \ref{lemma:Depen1}}  \label{subsubsec:lemmaDepen1} 
Applying \cite[Theorem 1.3]{Bosq}  with $b = 2$, $q = n^{(d-r)/(d + 1.25)}$,   and   $\eps =  \sqrt{32 r \log n}/\sqrt{q}$  gives 
$P\{|\mathrm{Re}(\varphi_n(t) - \varphi(t))|    \geq   \eps \}  \leq   \bar{o}(n^{-r})$ and  $P\{|\mathrm{Im}(\varphi_n(t) - \varphi(t))|    \geq   \eps \}  \leq   \bar{o}(n^{-r})$,   it thus follows 
\begin{equation}  \label{EqlemmaDepenAdd1}
P\{|\varphi_n(t) - \varphi(t)|    \geq   \sqrt{2} \eps \}  \leq  \bar{o}(n^{-r}). 
\end{equation} 

The remaining part of the proof is  similar to that of Lemma \ref{lemma:phi} so we keep it brief.  
Fix $\delta \in (1/2, \infty)$,  with the same  grid and similar arguments as in  Lemma \ref{lemma:phi},    it follows  that 
\begin{equation} \label{EqlemmaDepenAdd5}
P\{\sup_{\{ 0 \leq t  \leq    \log n\}}  |\varphi_n(t_k)  - \varphi(t)|  \geq (\sqrt{2}\eps  + \frac{1}{\sqrt{n}}) \} \leq I + II, 
\end{equation}  
where  $I  = P\{\sup_{\{1 \leq k  \leq n^{\delta} \log n\}}  |\varphi_n(t_k)  - \varphi(t_k)|  \geq   \sqrt{2} \eps \}$ and $II  \leq  P\{n^{-\delta} \sup_t  \{|\varphi'_n(t) - \varphi'(t)|\}   \geq \frac{1}{\sqrt{n}}\}$.  
The key for the proof is to show that 
\begin{equation}  \label{EqlemmaDepen2}
\mathrm{Var}(\frac{1}{n} \sum_{j = 1}^n |X_j|)  \leq C \log^{2a}(n)/n.   
\end{equation} 
In fact, once (\ref{EqlemmaDepen2}) is proved,  then on one hand,   by (\ref{EqlemmaDepenAdd1}), 
$I   \leq   n^{\delta} \log(n)   \cdot  \bar{o}(n^{-r})  =  \bar{o}(n^{-(r - \delta)})$.  On the other hand,        by  similar arguments as in the proof of Lemma \ref{lemma:phi},    
\[
II    \leq P\{\frac{1}{n}\sum_{j=1}^n  (|X_j| - E|X_j|) \geq n^{\delta - 1/2}  - s_n\}   \lesssim        
\frac{1}{n^{2\delta -1}}   \cdot   \mathrm{Var}(\frac{1}{n}\sum_{j=1}^n |X_j|) = \bar{o}(n^{-2 \delta}),      
\]
where $s_n \equiv \frac{1}{n} \sum_{j=1}^n E|X_j| \leq C \log^a(n)$  as  $\mmax_{j}\{|\mu_j| + |\sigma_j|\} \leq \log^a(n)$.    
The claims follows  by  taking $\delta = r/3$.  
  
We now show (\ref{EqlemmaDepen2}).    Applying  \cite[Corollary 1.1]{Bosq} with $p = 1.5$, $q = r = 6$,  
\[
\mathrm{Var}(\frac{1}{n} \sum_{j=1}^n |X_j|)    = \frac{1}{n^2}  \sum_{j, k}  \mathrm{Cov}(|X_j|, |X_k|)  \leq \frac{C}{n^2} \sum_{j, k}    \alpha^{2/3}(|j-k|)  \|X_j\|_6 \|X_k\|_6.   
\]
By  $\mmax_{j}\{|\mu_j| + |\sigma_j|\}  \leq  \log^a(n)$,      $\|X_j\|_6 \leq   C  \log^{a}(n)$ for all $1 \leq j \leq n$;  since $\alpha(k) \leq B k^{-d}$ with $d > 1.5$,     (\ref{EqlemmaDepen2}) follows by observing 
$\sum_{j, k}  \alpha^{2/3}(|j-k|) \leq C n   \sum_{k=1}^{\infty} \alpha^{2/3}(k)   \leq  C  n   \sum_{k=1}^{\infty}  k^{-2d/3}  = O(n)$.     \qed

\subsubsection{Proof of Lemma \ref{lemma:Depen2}}  \label{subsubsec:lemmaDepen2}
Consider the first claim.   By direct calculations,   
\[
|\varphi_n^{(k)}(t) -  \varphi^{(k)}(t)|   = |\frac{1}{n} \sum_{j = 1}^n (i X_j)^k e^{it X_j} - E[\frac{1}{n} \sum_{j = 1}^n (i X_j)^k e^{it X_j}]| \leq \frac{1}{n} \sum_{j = 1}^n [|X_j|^k + E|X_j|^k],  
\]
where the right hand side does not depend on $t$.  Since $\mmax_{\{j\}} \{|\mu_j| + |\sigma_j|\}  \leq B \log^a(n)$, the claim follows directly from  
$E|X_j|^{k}  \leq C \cdot  (|\mu_j|^k + |\sigma_j|^k)  \leq C \cdot   \log^{a k}(n)$,    $\forall  \; 1 \leq j \leq n$,   
where $C = C(k)$ is a generic constant. 

Consider the  second claim.        Introduce an event 
$D_n  = \{\mmax_{j}  \{|X_j|\} \leq 3B \log^{a + 1/2}(n)\}$.   
By  $\mmax_{\{j\}} \{|\mu_j| + |\sigma_j|\}  \leq B \log^a(n)$ and direct calculations,  
\begin{equation} \label{EqlemmaDepen7}  
P\{D_n^c\}  \leq   \sum_{j} P\{|X_j| \geq 3 B \log^{a+1/2}(n)\}    \leq  2  n  \bar{\Phi}(3 \sqrt{\log n} -1) =  \bar{o}(n^{-1}), 
\end{equation} 
where $\bar{\Phi}$ is the survival function of $N(0,1)$.    To show the claim, it suffices to show  
\begin{equation} \label{EqlemmaDepenAdd}
E[(\varphi_n^{(k)}(t_n)  - \varphi^{(k)}(t_n)) \cdot 1_{\{D_n\}}]^2   =  O(\log^{(2a+1)k}(n)/n). 
\end{equation}  
 
Now, first,  observe that   $|x|^k \mathrm{exp}(-\frac{(x - \mu_j)^2}{2\sigma_j^2}) = o(1)$, where $o(1) \goto  0$ as $n  \goto  \infty$,   uniformly for all $|x| \geq 3 B \log^{a + 1/2}(n)$ and $(\mu_j, \sigma_j)$ satisfying 
$|\mu_j| + |\sigma_j| \leq B \log^a(n)$;  combining this with (\ref{EqlemmaDepen7}) gives
$|E (\varphi_n(t_n) \cdot 1_{\{D_n^c\}})|    \leq \frac{1}{n} \sum_{j =1}^n  E(|X_j|^k \cdot 1_{\{D_n^c\}})   = \bar{o}(n^{-1})$.   
Notice that   $E\varphi_n^{(k)}(t_n)= \varphi^{(k)}(t_n)$, 
we thus have 
\begin{equation} \label{EqlemmaDepen8}
E[(\varphi_n^{(k)}(t_n)  -  \varphi^{(k)}(t_n)) \cdot 1_{\{D_n\}}]   =   -  E[\varphi_n^{(k)}(t_n) \cdot 1_{\{D_n^c\}}]     =    \bar{o}(1/n).   
\end{equation}  
Second,    as $\mmax_{\{j\}} \{|X_j|\} \leq 3 B \log^{a + 1/2}(n)$ over $D_n$,   by Billingsley's inequality  \cite[Page 22]{Bosq},   
\begin{align*}
  \mathrm{Var}(\varphi_n^{(k)}(t_n) \cdot 1_{\{D_n\}}) 
&=   \frac{1}{n^2} \sum_{j_1, j_2}    \mathrm{Cov}((i X_{j_1})^k \cdot  e^{it X_{j_1}}  \cdot 1_{\{D_n\}}, (i X_{j_2})^k \cdot  e^{i t X_{j_2}} \cdot 1_{\{D_n\}})     \\
& \leq  \frac{C}{n^2}  B^2 \log^{(2a+1)k}(n)   \sum_{j_1, j_2}   \alpha(|j_1-j_2|)   \leq O(\log^{(2a+1)k}(n)/n).    
\end{align*} 
Combining  this with (\ref{EqlemmaDepen8})  gives (\ref{EqlemmaDepenAdd}).  \qed


\subsection{Proof of Lemma \ref{lemma:Est1}}
For short, we drop $t$ from the functions whenever there is no confusion.   For the first claim,   by direct calculations, we have: 
\[
\so^2(g,t) - \so^2(f,t)  =   \frac{\frac{d}{dt}|f|}{t|f|}  -   \frac{\frac{d}{dt}|g|}{t|g|} =   I +   II + III,
\]     
where   $I  = (1 - \frac{|g|^2}{|f|^2}) \cdot \so^2(g,t)$,  
$ II = \frac{1}{ t |f|^2} \cdot    [\mathrm{Re}(g') \cdot  \mathrm{Re}(f  - g)  +  \mathrm{Im}(g') \cdot  \mathrm{Im}(f-g)  +    \mathrm{Re}(g)  \cdot  \mathrm{Re}((f - g)')  +  \mathrm{Im}(g)  \cdot  \mathrm{Im}((f - g)')]$,   
and 
$III    =    \frac{1}{t |f|^2}   \cdot  [\mathrm{Re}(f - g)  \cdot   \mathrm{Re}((f - g)')    +  \mathrm{Im}(f - g)  \cdot   \mathrm{Im}((f - g)')]$.     
Now, firstly,   using triangle inequality,  
\[
|I| \leq   \frac{|\so^2(g,t)|}{|f|^2} \cdot  \big| |f|^2 - |g|^2 \big|  \leq   \frac{|\so^2(g,t)|}{|f|^2} \bigl( 2 |g| \cdot  |f - g| +  |f - g|^2 \bigr); 
\]  
secondly, using  Cauchy-Schwartz  inequality,   
$|\mathrm{Re}(z) \mathrm{Re}(w) + \mathrm{Im}(z) \mathrm{Im}(w)| \leq |z| \cdot |w|$  for any complex numbers $z$ and $w$,   so it follows that
\[
|II| \leq \frac{1}{t |f|^2}   \cdot [|g'| \cdot |f  - g| + |g| \cdot |(f  - g)'|],  \qquad 
|III|  \leq   \frac{1}{t |f|^2} \cdot  |f - g| \cdot  |(f - g)'|;
\]   
combining these gives 
\[
|\so^2(g,t) - \so^2(f,t)| \leq \frac{1}{t |f|^2}   \bigl[(2 t \cdot  |\so^2(g,t)| \cdot  |g| + |g'|)  |f -g| + |g| \cdot  |(f - g)'| + \tilde{r}_n^{(1)} \bigr],
\]    
where  
$\tilde{r}_n^{(1)} =  t   \cdot  |\so^2(g,t)| \cdot    |f -g|^2  +  |f - g| \cdot |(f- g)'|$,   
and the claim follows directly.  

For the second claim,   by direct calculations:
\begin{align*}
\uo(g,t) - \uo(f,t) &= \frac{\mathrm{Re}(f')  \mathrm{Im}(f)  -   \mathrm{Re}(f) \mathrm{Im}(f')}{|f|^2}   -  \frac{\mathrm{Re}(g')  \mathrm{Im}(g)  -   \mathrm{Re}(g) \mathrm{Im}(g')}{|g|^2} \\
& = I + II + III,
\end{align*}
where  $I  =   (1  - \frac{|g|^2}{|f|^2}) \cdot  \uo(g,t)$,      
$II  =  \frac{1}{|f|^2} \cdot [(\mathrm{Re}(g') \cdot  \mathrm{Im}(f -   g)  -   \mathrm{Im}(g')  \cdot \mathrm{Re}(f - g))   +  (\mathrm{Im}(g) \cdot \mathrm{Re}((f - g)') -  \mathrm{Re}(g) \cdot \mathrm{Im}((f -  g)'))]$, and  
$III  = \frac{1}{|f|^2}  [\mathrm{Re}((f  - g)')  \cdot  \mathrm{Im}(f- g)  -   \mathrm{Re}(f - g) \cdot \mathrm{Im}((f-  g)')]$.
As in the first part,  
\[
|I| \leq  \frac{|\uo(g,t)|}{|f|^2}  [2   |g| \cdot  |f - g|  + |f - g|^2],
\]
\[
|II| \leq  \frac{1}{|f|^2} \cdot [ |g'| \cdot |f  -  g|  + |g| \cdot |(f  - g)'|],  \qquad   
|III| \leq   \frac{1}{|f|^2} \cdot |(f - g)'| \cdot |f  - g|;
\]   
combining these gives 
\[
|\uo(g, t) - \uo(f,t)|  \leq  \frac{1}{|f|^2} \cdot \bigl[(2 |\uo(g,t)| \cdot  |g| + |g'|) \cdot |f - g|  +   |g| \cdot |(f - g)'|  +  \tilde{r}_n^{(2)} \bigr],
\]   
where  
$\tilde{r}_n^{(2)}  =    |\uo(g,t)|  \cdot  |f - g|^2 + |f - g| \cdot |(f -g)'|$,      
and the claim follows. 
\qed

\subsection{Proof of Lemma \ref{lemma:phi}}
Lay out a grid $t_k  = k/n^{\delta}$,   for $k  = 1, \ldots, n^{\delta} \log n$ and $\delta \in (1/2, q_1/2)$.      For any $0 \leq t \leq  \log n$, pick the closest grid point $t_k$, so that $|t_k  - t| \leq n^{-\delta}$ and  
\[
|\varphi_n(t) - \varphi(t)|  \leq   |\varphi_n(t_k) - \varphi(t_k)|  +   |(\varphi_n(t) - \varphi(t)) - (\varphi_n(t_k) - \varphi(t_k))|, 
\] 
where  the second term is 
$\leq   n^{-\delta} \cdot  \sup_{t}    |\varphi'_n(t) - \varphi'(t)|$.  Write:
\[
\frac{\sqrt{2q_1 \log n} }{\sqrt{n}}   =  \lambda_1(q_1,n)  + \lambda_2(q_1,n),  
\]
where 
$\lambda_1(q_1, n)  =       \frac{\sqrt{2q_1  \log n}   - 2\log  \log n/\sqrt{2q_1 \log n}}{\sqrt{n}}$ and    $\lambda_2(q_1, n)  =    \frac{2\log  \log n/\sqrt{2q_1 \log n}}{\sqrt{n}}$.  
It thus follows  that 
\begin{equation}   \label{EqConfPf1}
P\{ \sup_{ 0 \leq t \leq \log n}  |\varphi_n(t) - \varphi(t)|  \geq \frac{\sqrt{2q_1 \log n} }{\sqrt{n}}  \}  \leq   I + II, 
\end{equation} 
where 
$I =  P\{\sup_{1 \leq k  \leq  n^{\delta} \log n}   |\varphi_n(t_k) - \varphi(t_k)| \geq  \lambda_1(q_1,n)\}$, and   
$II  =  P\{ n^{-\delta}  \cdot   \sup_{t}   |\varphi'_n(t) - \varphi'(t)| \geq  \lambda_2(q_1,n)\}$.     

For I,   a direct generalization of  Hoeffding's  inequality \cite{Hoeffding} to complex-valued random variables gives:   
\begin{align}  
I   & \leq     (n^{\delta}  \log n)  4 e^{- \frac{1}{4}n \lambda_1^2(q_1,n)}  = 4 n^{\delta} \log n  \cdot  e^{- \frac{q_1 \log n}{2} +  \log \log n (1 - \frac{\log \log n}{2q_1 \log n})}     \label{EqConfPf2add1}  \\
&\lesssim   (4 n^{\delta} \log n) (n^{-q_1/2} \log n) = 4 n^{\delta - q_1/2} \log ^2n.     \label{EqConfPf2add2}   \end{align} 
For II,   direct calculations show that  
$\sup_{t} |\varphi'_n(t) - \varphi'(t)|  \leq   \frac{1}{n} \cdot \sum_{j = 1}^n  (|X_j| + E|X_j|)$.   
Denote  $s_n   = \frac{1}{n} \sum_{j = 1}^n E|X_j|$ for short,   it follows from Chebyshev's inequality that:  
\begin{align} 
II  &\leq  P\{\frac{1}{n} \sum_{j=1}^n  (|X_j| + E|X_j|)  \geq    n^{\delta} \cdot  \lambda_2(q,n)\}       \label{EqConfPf3add1}   \\
&=   P\{     \frac{1}{n}  \sum_{j = 1}^n   (|X_j| - E|X_j|) \geq    n^{\delta}  \cdot  \lambda_2(q,n)  - 2 s_n \}     =    O(n^{-2 \delta}  \frac{\log^2(\log (n))}{\log(n)}),      \label{EqConfPf3add2}  
\end{align} 
where we have used the fact that $s_n$ is uniformly bounded from above by a constant $C(q,A,\uo, \so)  < \infty$.  
Inserting (\ref{EqConfPf2add1}) - (\ref{EqConfPf3add2})   to (\ref{EqConfPf1}) and taking $\delta = q_1/6$ give: 
\[
P\{ \sup_{ 0 \leq t \leq \log n}  |\varphi_n(t) - \varphi(t)|  \geq \frac{\sqrt{ 2q_1 \cdot   \log n} }{\sqrt{n}}  \}  = 4  \log^2(n) \cdot n^{-   q_1/3}  \cdot  (1 + o(1)),   \qquad  q_1 > 3. 
\]
This concludes the proof of Lemma \ref{lemma:phi}. \qed

\subsection{Proof of Lemma \ref{lemma:t}.}  
For short, write $\that_n = \that_n(\gamma)$, $t_n = t_n(\gamma)$,   $\varphi_n(t) = \varphi_n(t; X_1, \ldots, X_n,n)$,   $\varphi(t) = \varphi(t; \mu, \sigma, n)$, and $\Lambda_n = \Lambda_n(q,A; \uo, \so, \eps_0)$.    We claim that for  sufficiently large $n$,   $|\varphi(t)|$ is monotonely decreasing in $t$ over  $[\log \log n, \infty)$.   In fact,  
using  Lemma \ref{lemma:g},     when $n \goto \infty$,     
$\inf_{\{t \geq  \log \log n \}} \{ \so^2(\varphi; t) \}  = \so^2 \cdot (1 + o(1))  > 0$; 
recall that 
\begin{equation}  \label{Eqmonotonicityadd}
 \frac{d}{dt}|\varphi(t)| =   - t \cdot  |\varphi(t)| \cdot   \so^2(\varphi,t),   
\end{equation}    
the monotonicity  follows directly.

We now focus on the event  $D_n  = \{W_0(\varphi_n; n) \leq \sqrt{2 q_1 \log n}/\sqrt{n}\}$.  Recall that  $|\varphi(t_n)| = |\varphi_n(\that_n)| = n^{-\gamma}$,  so
\begin{equation}  \label{EqlemmatPf}
 \big| |\varphi(\that_n)|  - |\varphi(t_n)| \big|  =   \big| |\varphi(\that_n)|  - |\varphi_n(\that_n)| \big|   \leq  |\varphi(\that_n)  - \varphi_n(\that_n)|  \leq \sqrt{2 q_1 \log n}/\sqrt{n};    
\end{equation}
combining (\ref{Eqmonotonicityadd}) and (\ref{EqlemmatPf}) and 
using Taylor expansion,  there is a $\xi$ falling between $t_n$ and $\that_n$ such that
\begin{equation}  \label{Eqlemmatpf2}
|\that_n  - t_n|   =   \bigg|  \frac{|\varphi(\that_n)| - |\varphi(t_n)|}{|\varphi|'(\xi)}  \bigg|    \leq    \frac{\sqrt{2 q_1 \log n}/\sqrt{n}}{\xi \cdot  |\varphi(\xi)| \cdot  |\so^2(\varphi, \xi)|}.  
\end{equation}

At the same time,   elementary calculus shows 
\begin{equation}  \label{EqlemmaPf1}
(1 - 2 \eps_0)  e^{-\so^2 t^2 /2}  \leq  |\varphi(t)|  \leq e^{-\so^2 t^2/2},       \qquad   \forall \;  t > 0. 
\end{equation}  
Combining (\ref{EqlemmatPf}) and    (\ref{EqlemmaPf1}), it follows that     $\that_n \geq \log \log n$    for sufficiently large $n$.     Since $|\varphi(t)|$ is monotone over $[\log \log n, \infty)$,   so (\ref{EqlemmatPf}) and (\ref{EqlemmaPf1}) further imply that 
$ |\varphi(\xi)| \sim n^{-\gamma}$ and $\xi \sim  \that_n  \sim t_n \sim \sqrt{2 \gamma \log n}/\so$;  these, together with Lemma \ref{lemma:g},   imply that 
$\so^2(\varphi, \xi) \sim \so^2$.     
Inserting these into (\ref{Eqlemmatpf2}) gives
$|\that_n  - t_n|   \lesssim    \frac{\sqrt{2 q_1 \log n}/\sqrt{n}}{\so^2 \cdot  t_n \cdot n^{\gamma}} \sim   \frac{1}{\so} \cdot \sqrt{q_1/\gamma} \cdot n^{\gamma -1/2}$.       \qed

\subsection{Proof of Lemma \ref{lemma:phi'}}
Lay out a grid  $t_k = (t_n(\gamma)  -  \tau_0  n^{\gamma -1/2})  + \frac{k}{n^{\delta}}$,    for $1 \leq k \leq  2  \tau_0   n^{\delta + \gamma - 1/2}$ and $\delta \in [1/2, \infty)$.       For any $t \in [t_k,  t_{k + 1}]$,  \begin{equation}
|\varphi_n'(t)  - \varphi'(t)|  \leq   |\varphi_n'(t_k)  - \varphi'(t_k)|  +   n^{-\delta} \cdot \biggl( \sup_{|\xi - t_n(\gamma)| \leq \tau_0 \cdot n^{\gamma  - 1/2} }|\varphi_n''(\xi) - \varphi''(\xi)| \biggr).      
\end{equation}  
By direct calculations and the definition of $\bar{s}_n$,   
\[
|\varphi''_n(\xi) - \varphi''(\xi)|  \leq  \frac{1}{n} \sum_{j = 1}^n  (X_j^2 + E[X_j^2])  \equiv  \frac{1}{n} \sum_{j = 1}^n  (X_j^2 - E[X_j^2])    + 2 \bar{s}_n^2,   
\]
it thus follows that:  
\begin{align}    
|\varphi_n'(t)  - \varphi'(t)|   & \leq       |\varphi_n'(t_k)  - \varphi'(t_k)|    + n^{-\delta} \cdot  \biggl[\frac{1}{n} \sum_{j = 1}^n (X_j^2   -  E(X_j^2))  + 2 \bar{s}_n^2\biggl]    \label{Eqphi'Add11}  \\
&\leq   |\varphi_n'(t_k)  - \varphi'(t_k)|     +   n^{-\delta} \cdot [\frac{1}{n} \sum_{j = 1}^n (X_j^2   -  E(X_j^2))]  + \frac{2 \bar{s}_n^2}{\sqrt{n}}.    \label{Eqphi'Add12} \end{align} 

Now,  denote $q_1 = q/2 -1$ for short,   write:
\begin{equation}   \label{Eqphi'Add2}
\frac{\bar{s}_n    (\sqrt{(q-2) \log n}  +  2 \bar{s}_n)}{\sqrt{n}}   = \frac{\bar{s}_n    (\sqrt{2 q_1 \log n}  +  2 \bar{s}_n)}{\sqrt{n}}  =   \lambda_1(q,n)  + \lambda_2(q,n) + \frac{2 \bar{s}_n^2}{\sqrt{n}},   
\end{equation} 
where  $\lambda_1(q_1,n) = \bigl(\bar{s}_n  \sqrt{2 q_1 \log n}   -  (\frac{\log \log n}{2 \bar{s}_n \sqrt{2q_1  \log n}})\bigr)/\sqrt{n}$ and  $\lambda_2(q_1,n)  =  \bigl(\frac{ \log \log n}{2 \bar{s}_n \sqrt{2q_1  \log n}} \bigr)/\sqrt{n}$. 
Compare  (\ref{Eqphi'Add2}) with  (\ref{Eqphi'Add11}) - (\ref{Eqphi'Add12})  gives: 
\[
P\{\sup_{|t - t_n(\gamma)| \leq  \pi_0  \cdot n^{\gamma - 1/2}}   |\varphi_n'(t)  - \varphi'(t)|  \geq \frac{\bar{s}_n \cdot (\sqrt{2 q_1 \log n} + 2 \bar{s}_n)}{\sqrt{n}} \}   \leq I + II, 
\]
where 
$I   = P\{ \sup_{1 \leq k \leq 2  \pi_0  n^{\delta + \gamma - 1/2}}    |\varphi_n'(t_k)  - \varphi'(t_k)|  \geq   \lambda_1(q_1,n)\}$,      
and 
$II  =   P\{n^{-\delta} \cdot [\frac{1}{n} \sum_{j = 1}^n  (X_j^2  - E X_j^2)] \geq   \lambda_2(q_1,n) \}$.      

For I,  by   \cite[Theorem 1]{Rubin} and direct calculations,   
\begin{equation}  \label{Eqphi'pf3add1}
I \leq   (2 \pi_0 n^{\delta + \gamma -1/2})  \cdot \bar{o}(e^{- \frac{1}{2}n \lambda_1^2(q_1,n)})    \leq (2 \pi_0 n^{\delta + \gamma -1/2})  \cdot  \bar{o}(n^{-q_1})  = \bar{o}(n^{\delta + \gamma - 1/2 - q_1}).      
\end{equation}
For II, we study for the case $q < 4$ and the case $q \geq 4$ separately.  

For  the case $q < 4$,     set  $\delta = (q_1 + 1- \gamma)/2 > 1/2$,   by   Chebyshev's  inequality, 
\begin{equation}   \label{Eqphi'pf3}
II  =   P\{\frac{1}{n} \sum_{i = 1}^n  X_j^2    \geq   \bar{s}_n^2 +  n^{\delta} \cdot \lambda_2(q_1,n) \} \leq    \frac{ \bar{s}_n^2 }{\bar{s}_n^2  + n^{\delta} \cdot  \lambda_2(q_1,n)}   \leq  \bar{o}(n^{1/2 - \delta}),   
\end{equation} 
where we have used the fact that  $\bar{s}_n^2$ is   uniformly bounded from  above by a constant $C_1 = C_1(q,A, \uo, \so) < \infty$.  
Notice that the choice of $\delta$ satisfies  
$\delta + \gamma - 1/2 -  q_1 = 1/2 - \delta = (1 + \gamma  - q/2)/2$,       
combining (\ref{Eqphi'pf3add1}) and (\ref{Eqphi'pf3})   gives   $I + II \leq \bar{o}(n^{(1 +  \gamma - q/2)/2})$.

For the case $q \geq 4$,    notice that  $\frac{1}{n}  \sum_{j = 1}^n  E(X_j^2 - E[X_j^2])^2$ is  uniformly bounded from above by a constant $C_2 = C_2(q,A,\uo,\so) < \infty$,       it  follows from  Chebyshev's inequality that
\begin{equation}  \label{Eqphi'pf3add2} 
 II \leq \biggl(\frac{C_2}{\lambda_2^2(q_1,n) \cdot n \cdot n^{2 \delta}}\biggr)  =    \bar{o}(n^{-2 \delta}).  
\end{equation} 
Set $\delta  =   \mmax\{1/2,  (q - 1 - 2 \gamma)/6\}$,     combining (\ref{Eqphi'pf3add1}) and (\ref{Eqphi'pf3add2}) gives:       
\[
I  + II \leq  \left\{
\begin{array}{ll}
 \bar{o}(n^{ \gamma + 1 - q/2}),  &\qquad  4 \leq q  \leq 4 + 2 \gamma,  \\
 \bar{o}(n^{(2 \gamma + 1 - q)/3}),  &\qquad  q > 4 + 2 \gamma.   
 \end{array}
 \right.
\]
This  finishes the proof of Lemma \ref{lemma:phi'}.      \qed
\subsection{Proof of Lemma \ref{lemma:g}}
First, we show   (\ref{lemmag1}).  Write    $|\varphi(t)| = |\varphi_0(t)| \cdot |1 + r(t)|$,  recall   that  $\so^2(\varphi_0; t) \equiv \so^2 t$,  so     
\[
\frac{\frac{d}{dt} |\varphi(t)|}{|\varphi(t)|}  =  \frac{\frac{d}{dt}|\varphi_0(t)| \cdot |1 + r(t)| + |\varphi_0(t)| \cdot  \frac{d}{dt}|1 + r(t)|}{|\varphi_0(t)| \cdot |1 + r(t)|}    =     -\so^2 t  + \frac{\frac{d}{dt}|1 + r(t)|}{|1 + r(t)|}, 
\]
and  it follows that 
$\so^2(\varphi,t) - \so^2 =  -  \frac{d}{dt}(|1 + r(t)|)/(t \cdot |1 + r(t)|)$,   
which yields  (\ref{lemmag1})   by direct calculations.

Next, we show (\ref{lemmag2}).   For short,   we drop $t$ from all  expressions whenever there is no confusion.    Since    $\varphi = \varphi_0(1 + r)$,        
 $\mathrm{Re}(\varphi) =   \mathrm{Re}(\varphi_0)   +   \mathrm{Re}(r)  \mathrm{Re}(\varphi_0)   - \mathrm{Im}(r) \mathrm{Im}(\varphi_0)$, and    $\mathrm{Im}(\varphi) =  \mathrm{Im}(\varphi_0)  +  \mathrm{Im}(r)  \mathrm{Re}(\varphi_0)    +  \mathrm{Re}(r)   \mathrm{Im}(\varphi_0)$;   it can be showed that 
\begin{equation}  \label{EqApprox1aadd1}
\mathrm{Re}(\varphi')   \cdot    \mathrm{Im}(\varphi)  -  \mathrm{Im}(\varphi')   \cdot   \mathrm{Re}(\varphi)   =  I + II,  
\end{equation} 
where   
$I  =    - |1 +  r|^2 \uo |\varphi_0|^2$,    and  
$II =  |\varphi_0|^2 \cdot [ -  \mathrm{Im}(r')  +   \mathrm{Re}(r') \mathrm{Im}(r)  -  \mathrm{Im}(r') \mathrm{Re}(r)]$.  The proof of (\ref{EqApprox1aadd1}) is long, so we leave it to the end of this section.  Now,  
\begin{align*} 
\uo(\varphi;t)   =  - \frac{I  + II}{|\varphi|^2}   =\uo +   \frac{\mathrm{Im}(r')  -    \mathrm{Re}(r') \mathrm{Im}(r)  + \mathrm{Re}(r) \mathrm{Im}(r')}{|1 + r|^2},      
\end{align*}
so by  Cauchy-Schwartz  inequality,   
\[
|\uo(\varphi,t) - \uo|  =   \frac{ | \mathrm{Im}(r')  -    \mathrm{Re}(r') \mathrm{Im}(r)  +  \mathrm{Im}(r') \mathrm{Re}(r) | }{|1 +  r|^2} \leq  |r'| \cdot \frac{1 + |r|}{|1+r|^2},   
\]
and (\ref{lemmag2})  follows directly.    

Next, we show (a1) and (a3). (a2) follows directly from (a1) and direct calculations, so we omit it.  

We now show (a1).  For the  $5$ inequalities, the proofs for the first $4$ are similar, so we only show the second one and the last one.    First,  consider  the second inequality.    Use  H\"older's inequality,   
\begin{align} 
\mathrm{Ave}_{\{j: \; (\mu_j, \sigma_j) \neq (\uo, \so)\}}   \{|\mu_j - \uo|  + (\sigma_j^2 - \so^2)^{1/2}\}    &\leq A,    \label{Eqlemmagadd21} \\
\mathrm{Ave}_{\{j: \; (\mu_j, \sigma_j) \neq (\uo, \so)\}}   \{(\sigma_j^2 - \so^2)\}   &\leq A^2.     \label{Eqlemmagadd22}
\end{align} 
Note that  $\sup_{\{x \geq 0\}}  x e^{-x^2/2} = 1/e \leq 1$,  direct calculations show that
\begin{align}  
|g'(t)|  &\leq  \mathrm{Ave}_{\{j: \; (\mu_j, \sigma_j) \neq (\uo, \so)\}}   \{ e^{ - \frac{(\sigma_j^2 - \so^2)t^2}{2}}  \cdot  [|\mu_j - \uo|  + (\sigma_j^2 - \so^2) t] \} 
\label{g1}  \\
&\leq  \mathrm{Ave}_{\{j: \; (\mu_j, \sigma_j) \neq (\uo, \so)\}} \{|\mu_j - \uo| +  (\s_j^2 - \so^2)^{1/2}  \cdot   [(\s_j^2 - \so^2)^{1/2} t \cdot e^{-\frac{(\s_j^2 - \so^2) t^2}{2}}]\}    \nonumber   \\
&   \leq   \mathrm{Ave}_{\{j: \; (\mu_j, \sigma_j) \neq (\uo, \so)\}} \{ |\mu_j - \uo|  + (\s_j^2 - \so^2)^{1/2}\},       \nonumber 
\end{align}  
 the second inequality  follows directly by using  (\ref{Eqlemmagadd21}).     
Second, consider the last inequality.  
By the definition of $\tau_n$ and (\ref{Eqlemmagadd21}),  it is seen 
\begin{equation}   \label{Eqlemmagadd31} 
\mathrm{Ave}_{\{j: \; (\mu_j, \sigma_j) \neq (\uo, \so)\}}   \{ e^{ - \frac{(\sigma_j^2 - \so^2)t^2}{2}}  \cdot  |\mu_j - \uo| \}     
 \leq  A e^{-\frac{\tau_n t^2}{2}}.   
\end{equation} 
At the same time,   notice that $\sup_{\{x \geq 0\}}  x e^{-x/2} = 2/e$,   so 
$e^{ - \frac{(\sigma_j^2 - \so^2)t^2}{2}}  \cdot  (\sigma_j^2 - \so^2) t \leq  
\mmin\{e^{-\tau_n t^2/2}  \cdot  (\sigma_j^2 - \so^2) t ,   2/(et)\}$,    and it follows from     (\ref{Eqlemmagadd22}) that  
\begin{equation}   \label{Eqlemmagadd4}   
\mathrm{Ave}_{\{j: \; (\mu_j, \sigma_j) \neq (\uo, \so)\}}\{e^{ - \frac{(\sigma_j^2 - \so^2)t^2}{2}}  \cdot  (\sigma_j^2 - \so^2) t\} \leq    \mmin\{A^2 e^{-\tau_n t^2/2} t,  2/(et)\}.   
\end{equation} 
The claim follows by combining (\ref{g1}) -   (\ref{Eqlemmagadd4}).

Next, we show (a3).   As the proofs are similar, we only show that corresponds to $\so^2$.   By  (\ref{lemmag1}),  $|\so^2(\varphi; t) - \so^2| \goto 0$ uniformly;  by (a1),  it is not hard to show that $\so^2(\varphi; t)$ and its first two  derivatives are all  uniformly bounded;    so  all remains to show is that $\frac{d}{dt}\so^2(\varphi; t) \goto 0$ uniformly.     Observe that for any twice differentiable function $f$ and  $\Delta > 0$,      $|\frac{f(t + \Delta) - f(t)}{\Delta} - f'(t)| \leq  \sup_{\{s\}} \{|f''(s)|\} \Delta$,   so it follows $|f'(t)| \leq    \{ \sup_{\{s\}} \{|f''(s)|\} \Delta   +  \frac{1}{\Delta} \sup_{\{s,s' \geq t\}}\{ |f(s) - f(s')|\}$;   the claim follows by  taking $\Delta = \sqrt{\sup_{\{s,s' \geq t\}}\{ |f(s) - f(s')|\}/\sup_{\{s\}} \{|f''(s)|\}}$ and $f(t) = \so^2(\varphi; t)$.  
 
Lastly, we validate  (\ref{EqApprox1aadd1}).        
Write $\mathrm{Re}(\varphi')   =  \mathrm{Re}(\varphi_0')   +   \mathrm{Re}(r')  \mathrm{Re}(\varphi_0)     +  \mathrm{Re}(r)  \mathrm{Re}(\varphi_0')  - \mathrm{Im}(r') \mathrm{Im}(\varphi_0) - \mathrm{Im}(r) \mathrm{Im}(\varphi_0')$, and $\mathrm{Im}(\varphi')   =   \mathrm{Im}(\varphi_0')  +  \mathrm{Im}(r')  \mathrm{Re}(\varphi_0)  +  \mathrm{Im}(r)  \mathrm{Re}(\varphi_0')    +  \mathrm{Re}(r')   \mathrm{Im}(\varphi_0) + \mathrm{Re}(r)   \mathrm{Im}(\varphi_0')$, 
we have
\begin{align*}
& \mathrm{Re}(\varphi')   \cdot    \mathrm{Im}(\varphi)  -  \mathrm{Im}(\varphi')   \cdot   \mathrm{Re}(\varphi)       =     \mathrm{Re}(\varphi_0')  [ \mathrm{Im}(\varphi_0)  +  \mathrm{Im}(r)  \mathrm{Re}(\varphi_0)    +  \mathrm{Re}(r)   \mathrm{Im}(\varphi_0)]   \\
& + \mathrm{Re}(r') \mathrm{Re}(\varphi_0)   [\mathrm{Im}(\varphi_0)  +  \mathrm{Im}(r)  \mathrm{Re}(\varphi_0)  +  \mathrm{Re}(r)   \mathrm{Im}(\varphi_0)]   
+  \mathrm{Re}(r)  \mathrm{Re}(\varphi_0') [\mathrm{Im}(\varphi_0) 
+  \\ &    \mathrm{Im}(r)    
\mathrm{Re}(\varphi_0)     +  \mathrm{Re}(r)   \mathrm{Im}(\varphi_0)]  
 -  \mathrm{Im}(r') \mathrm{Im}(\varphi_0)    [\mathrm{Im}(\varphi_0)   
+  \mathrm{Im}(r)  \mathrm{Re}(\varphi_0)   +  \mathrm{Re}(r)   
\mathrm{Im}(\varphi_0)]  \\&     -   \mathrm{Im}(r) \mathrm{Im}(\varphi_0')   
  [\mathrm{Im}(\varphi_0)  +  \mathrm{Im}(r)  \mathrm{Re}(\varphi_0)   
+  \mathrm{Re}(r)   \mathrm{Im}(\varphi_0)]  -  \mathrm{Im}(\varphi_0')  [\mathrm{Re}(\varphi_0)    +   \mathrm{Re}(r)   \\&  \mathrm{Re}(\varphi_0)   
   - \mathrm{Im}(r) \mathrm{Im}(\varphi_0)]   
  -  \mathrm{Im}(r') \mathrm{Re}(\varphi_0)   [\mathrm{Re}(\varphi_0)   +   \mathrm{Re}(r)  \mathrm{Re}(\varphi_0)   - \mathrm{Im}(r)   
   \mathrm{Im}(\varphi_0)]  -    \\&  
 \mathrm{Im}(r) \mathrm{Re}(\varphi_0')  
  [\mathrm{Re}(\varphi_0)    
  +   \mathrm{Re}(r)  \mathrm{Re}(\varphi_0)   - \mathrm{Im}(r) \mathrm{Im}(\varphi_0)] - \mathrm{Re}(r') \mathrm{Im}(\varphi_0)   
[\mathrm{Re}(\varphi_0)   +   
\\ &  \mathrm{Re}(r)  \mathrm{Re}(\varphi_0)  
 -  \mathrm{Im}(r) \mathrm{Im}(\varphi_0)]   
  -   \mathrm{Re}(r) \mathrm{Im}(\varphi_0')   [\mathrm{Re}(\varphi_0)   +   \mathrm{Re}(r)  \mathrm{Re}(\varphi_0)
   - \mathrm{Im}(r) \mathrm{Im}(\varphi_0)]; 
\end{align*}
by cancellations,  this reduces to
 \begin{align*}
&  \mathrm{Re}(\varphi_0') \cdot [ \mathrm{Im}(\varphi_0)   +  \mathrm{Re}(r)   \mathrm{Im}(\varphi_0)]  
+ \mathrm{Re}(r') \mathrm{Re}(\varphi_0)  [ \mathrm{Im}(r)  \mathrm{Re}(\varphi_0)]   
+  \mathrm{Re}(r)  \mathrm{Re}(\varphi_0') \\
&    [\mathrm{Im}(\varphi_0)     +  \mathrm{Re}(r)   \mathrm{Im}(\varphi_0)]  -  \mathrm{Im}(r') \mathrm{Im}(\varphi_0)     [\mathrm{Im}(\varphi_0)    +  \mathrm{Re}(r)   \mathrm{Im}(\varphi_0)] 
 -   \mathrm{Im}(r) \mathrm{Im}(\varphi_0')   \\
 & [ \mathrm{Im}(r)  \mathrm{Re}(\varphi_0)]   
-  \mathrm{Im}(\varphi_0') \cdot [\mathrm{Re}(\varphi_0)   +   \mathrm{Re}(r)  \mathrm{Re}(\varphi_0)] 
-  \mathrm{Im}(r') \mathrm{Re}(\varphi_0)    [\mathrm{Re}(\varphi_0)  
  +   \mathrm{Re}(r)   \\
& \mathrm{Re}(\varphi_0)]  -  \mathrm{Im}(r) \mathrm{Re}(\varphi_0') \cdot [   - \mathrm{Im}(r) \mathrm{Im}(\varphi_0)] 
- \mathrm{Re}(r') \mathrm{Im}(\varphi_0)    [ - \mathrm{Im}(r) \mathrm{Im}(\varphi_0)]  \\
&- \mathrm{Re}(r) \mathrm{Im}(\varphi_0')     [\mathrm{Re}(\varphi_0)   +   \mathrm{Re}(r)  \mathrm{Re}(\varphi_0)];   
\end{align*}
by recombinations,  this reduces to 
$|1 + r|^2 \cdot [\mathrm{Re}(\varphi_0')  \mathrm{Im}(\varphi_0)   
-  \mathrm{Re}(\varphi_0) \mathrm{Im}(\varphi_0')]
+  |\varphi_0|^2  \cdot [- \mathrm{Im}(r')  +   
\mathrm{Re}(r') \mathrm{Im}(r)   
   - \mathrm{Im}(r') \mathrm{Re}(r)] $.    
Note that  $[\mathrm{Re}(\varphi_0')  \mathrm{Im}(\varphi_0)   -  \mathrm{Re}(\varphi_0) \mathrm{Im}(\varphi_0')] = -\uo |\varphi_0|^2$,        (\ref{EqApprox1aadd1}) follows directly.   
 \qed

\end{document}